\documentclass[a4paper]{amsart}
\usepackage[T1]{fontenc}
\usepackage[utf8]{inputenc}
\usepackage[british]{babel}
\usepackage{amsfonts,amsmath,amssymb,amsthm,mathtools,esint}
\usepackage[dvipsnames]{xcolor}
\usepackage{tikz}
\usepackage{pgfplots}
\pgfplotsset{compat=1.3}
\usepackage[justification=centering,width=0.9\linewidth]{caption}
\usepackage{subcaption}
\usepackage{enumitem}
\usepackage[hidelinks]{hyperref}
\usepackage{cleveref}
\crefname{subsection}{subsection}{subsections}
\numberwithin{equation}{section}

\usepackage[backend=bibtex,style=numeric,natbib=true,doi=false,isbn=false,url=false,eprint=false,firstinits=true]{biblatex}

\newcommand{\D}{\mathrm{D}}

\newcommand{\vertiii}[1]
{{\left\vert\kern-0.25ex\left\vert\kern-0.25ex\left\vert #1 
\right\vert\kern-0.25ex\right\vert\kern-0.25ex\right\vert}}

\newcounter{cnst}
\newcommand{\newcnst}{%
	\refstepcounter{cnst}%
	\ensuremath{C_{\thecnst}}}
\newcommand{\cnst}[1]{\ensuremath{C_{\ref{#1}}}}

\newtheorem{theorem}{Theorem}[section]
\newtheorem{proposition}[theorem]{Proposition}
\newtheorem{lemma}[theorem]{Lemma}

\newtheorem{assumption}[theorem]{Assumption}
\newtheorem{corollary}[theorem]{Corollary}
\theoremstyle{remark}
\newtheorem{remark}[theorem]{Remark}

\begin{document}

\title[FEM for linear PDE in nondivergence form]{Finite element approximation for uniformly elliptic linear PDE of second order in nondivergence form}

\author[N.~T.~Tran]
{Ngoc Tien Tran}
\thanks{This project received funding from
	the European Union's Horizon 2020 research and innovation 
	programme (project DAFNE, grant agreement No.~891734, and project RandomMultiScales, grant agreement No.~865751).
	}
\address[N.~T.~Tran]{%
	Institut f\"ur Mathematik,
	Uni\-ver\-si\-t\"at Augsburg,
	Universit\"atsstra\ss e 2, 86159 Augsburg, Germany}
\email{ngoc1.tran@uni-a.de}
\date{\today}



\begin{abstract}
	This paper proposes a novel technique for the approximation of strong solutions $u \in C(\overline{\Omega}) \cap W^{2,n}_\mathrm{loc}(\Omega)$ to uniformly elliptic linear PDE of second order in nondivergence form with continuous leading coefficient in nonsmooth domains by finite element methods (FEM).
	These solutions satisfy the Alexandrov-Bakelman-Pucci (ABP) maximum principle, which provides an a~posteriori error control for $C^1$ conforming approximations.
	By minimizing this residual,
	we obtain an approximation to the solution $u$ in the $L^\infty$ norm.
	Although discontinuous functions do not satisfy the ABP maximum principle,
	this approach extends to nonconforming FEM as well thanks to well-established enrichment operators.
	Convergence of the proposed FEM is established for uniform mesh-refinements.
	The built-in a~posteriori error control (even for inexact solve) can be utilized in adaptive computations for the approximation of singular solutions, which performs superiorly in the numerical benchmarks in comparison to the uniform mesh-refining algorithm.
\end{abstract}

\maketitle

\section{Introduction}
\subsection{Background}
Given an open bounded Lipschitz domain $\Omega \subset \mathbb{R}^n$, we seek the strong solution $u \in C(\overline{\Omega}) \cap W^{2,n}_\mathrm{loc}(\Omega)$ to the Dirichlet problem
\begin{align}\label{def:PDE}
	L u \coloneqq - A:\D^2 u + b \cdot \nabla u + c\,u = f \text{ in } \Omega \quad\text{and}\quad u = g \text{ on } \partial \Omega
\end{align}
with an uniformly elliptic second-order operator $L$ in nondivergence form, a right-hand side $f \in L^n(\Omega)$, and Dirichlet data $g \in C(\partial \Omega)$.
The existence of strong solutions is guaranteed under the following assumption.
\begin{assumption}\label{assumption:coefficients}
	Let $A \in C(\overline{\Omega};\mathbb{S})$ with $\lambda \mathrm{I}_n \leq A \leq \Lambda \mathrm{I}_n$, $b \in L^\infty(\Omega;\mathbb{R}^n)$, and $0 \leq c \in L^\infty(\Omega)$.
\end{assumption}
Here, $0 < \lambda \leq \Lambda$ are (fixed) ellipticity constants.
We refer to \Cref{sec:mathematical-setting} for further details on the PDE \eqref{def:PDE}.
Monotone finite difference methods (FDM) can approximate these solutions because they respect some maximum principle on the discrete level. General convergence theory has been established in \cite{BarlesSouganidis1991} even for fully nonlinear PDE.
However, FDM are restricted to low-order methods and
fixed size finite difference stencils are generally not sufficient for a convergent scheme \cite{MotzkinWasow1953}.
A remedy are wide stencil FDM, which are constructed on unstructured meshes in, e.g., \cite{DebrabantJakonsen2013,FengJensen2017} with convergence for fully nonlinear problems.

In contrast to PDE in divergence form, a variational formulation for \eqref{def:PDE} may not be available if $A$ is not sufficiently smooth. In these cases, the design of finite element schemes for \eqref{def:PDE} becomes challenging.
We point out several finite element methods (FEM) in the literature.
By imitating the convergence analysis on the continuous level, a FEM have been proposed in \cite{FengHenningsNeilan2017}.
The two-scale lowest-order method in \cite{NochettoZhang2018} relies on discrete maximum principles and allows for convergence of FEM under further assumptions on the mesh.
Both previously mentioned papers require the regularity $u \in W^{2,p}(\Omega)$ and therefore, $C^{1,1}$ boundary for the domain $\Omega$ as a sufficient condition.
For uniformly elliptic operators $L$ with (possibly) discontinuous coefficients that satisfy the so-called Cordes condition, \cite{SmearsSueli2013,SmearsSueli2014} prove that there exists a strong solution $u \in H^2(\Omega) \cap H^1_0(\Omega)$ to \eqref{def:PDE} for all right-hand side $f \in L^2(\Omega)$ on convex domains $\Omega$.
The PDE \eqref{def:PDE} can be seen as a perturbation of the Laplace equation as the eigenvalues of $A$ cannot spread too far depending on the dimension $n$.
This allows access to finite element discretization and adaptive computation with plain convergence in \cite{GallistlSueli2019,KaweckiSmears2022}.
Since the Cordes condition allows for discontinuous coefficients, which arise in the linearization of fully nonlinear Hamilton-Jacobi-Bellman and Isaacs equations, the aforementioned results also apply to these classes of fully nonlinear equations \cite{SmearsSueli2014,GallistlSueli2019,KaweckiSmears2021,KaweckiSmears2022}.
In two space dimensions and without lower-order terms, the Cordes condition therein is equivalent to the uniformly ellipticity of $L$. The restrictions imposed on the coefficients are less practical in higher space dimensions or in presence of lower-order terms.

Another approach related to this paper is the first order least-squares method from \cite{QiuZhang2020}. While it works well in the numerical examples presented therein, it raises several questions in regards of the theory.
Since this method mimics the design of least-squares schemes for PDE of second order in divergence form, the assumptions for the analysis therein are not tailored to the nondivergence case.
It is expected that these assumptions can only be verified for a very limited class of operators on nonsmooth domains.

\subsection{Motivation}\label{sec:motivation}
Global regularity of strong solutions $u$ to the PDE \eqref{def:PDE} can only be obtained under rather strict restrictions on the domain $\Omega$, e.g., $C^{1,1}$ boundary.
On Lipschitz domains, however, we can only expect the local regularity $u \in C(\overline{\Omega}) \cap W^{2,n}_\mathrm{loc}(\Omega)$; that is, $u \in W^{2,n}(\omega)$ for any open set $\omega \Subset \Omega$.
The goal of this paper is the design of convergent finite element methods on nonsmooth domains without additional assumptions.
The main tool is the well-known Alexandrov-Bakelman-Pucci (ABP) maximum principle
\begin{align}
	\|u\|_{L^\infty(\Omega)} &\leq \|u\|_V \coloneqq \|u\|_{L^\infty(\partial \Omega)} + \cnst{cnst:ABP}\|L u\|_{L^n(\Omega)}\label{ineq:ABP-L-infty}
\end{align}
for a positive constant $\newcnst\label{cnst:ABP}$ independent of $u$.
The key observation is that, given any function $v \in V$ in the Banach space
\begin{align}
	V \coloneqq \{v \in C(\overline{\Omega}) \cap W^{2,n}_\mathrm{loc}(\Omega) : L v \in L^n(\Omega)\}\label{def:V},
\end{align}
endowed with the norm $\|\bullet\|_V$ from \eqref{ineq:ABP-L-infty},
\eqref{ineq:ABP-L-infty} implies the stability estimate $\|u - v\|_{L^\infty(\Omega)} \leq \Phi(v)$ with the function
\begin{align}\label{def:Phi}
	\Phi(v) \coloneqq \|g - v\|_{L^\infty(\partial \Omega)} + \cnst{cnst:ABP}\|f - L v\|_{L^n(\Omega)},
\end{align}
which simultaneously provides an a~posteriori error control for the error of $u - v$ in the $L^\infty$ norm.
We will show that $W^{2,n}(\Omega)$ is dense in $V$ so that the infimum of $\Phi$ among all functions in $W^{2,n}(\Omega)$ vanishes.
In particular, the sequence of discrete minimizers $u_h$ of $\Phi$ in a $C^1$ conforming finite element space $V_h$ converges uniformly to $u$ as the mesh-size $h$ tends to zero.
Here, $u_h$ can be understood as a (possibly non-unique) best-approximation of $u$ in $V_h$ with respect to the norm $\|\bullet\|_V$.
The main difficulty of this fairly simple approach is the practical realization because $\Phi$ is a nonsmooth nonlinear functional.
If $g$ is the trace of a finite element function, it is possible to enforce the Dirichlet data pointwise onto $V_h$.
This leads to a smooth minimization problem in the affine space $W_h \coloneqq \{v_h \in V_h : v_h = g \text{ on } \partial \Omega\}$.
However, we will demonstrate with the Laplace equation as an example that this approach will fail in the sense that the minimum of $\Phi$ in $W_h$ may not vanish as $h \to 0$ and the sequence of minimizers of $\Phi$ in $W_h$ may not approximate $u$.
Instead, the boundary error $\|g - \bullet\|_{L^\infty(\partial \Omega)}$ in \eqref{def:Phi} is enforced as linear side constraints.
As a result, the proposed scheme requires solving a constrained convex minimization problem or, in two space dimensions $n = 2$, a quadratic programming.
While this is numerically less efficient than least-squares schemes, e.g., from \cite{QiuZhang2020}, convergence of this method is guaranteed for continuous coefficients $A$.

Due to their flexibility in terms of polynomial degree and simplicity of their practical realization, nonconforming discretizations may outperform conforming ones for problems involving the Hessian.
Although the ABP maximum principle \eqref{ineq:ABP-L-infty} cannot be directly applied to 
discontinuous functions,
an enrichment operator based on local averaging provides appropriate conforming approximations of these functions. Therefore, the convergence analysis of nonconforming FEM can be carried out as for conforming schemes.
A welcome feature of the analysis of this paper is the built-in a~posteriori error control that allows for adaptive mesh-refining strategies.

\subsection{Outline and Notation}
The remaining parts of this paper are organized as follows. \Cref{sec:mathematical-setting} recalls some classical results from PDE theory and proves the density of $W^{2,n}(\Omega)$ in $V$ with respect to the norm $\|\bullet\|_V$.
We demonstrate the design of FEM with this density result in \Cref{sec:convergencent-FEM} for conforming and nonconforming schemes.
Numerical benchmarks in \Cref{sec:numerical-examples} conclude this paper.

Standard notation for function spaces applies throughout this paper. Let $\mathbb{S} \subset \mathbb{R}^{n \times n}$ denote the set of all symmetric matrices with the identity matrix $\mathrm{I}_n$.
The notation $A : B$ denotes the Euclidean scalar product of two matrices $A, B \in \mathbb{R}^{n \times n}$, which induces the Frobenius norm $|\bullet|$ in $\mathbb{R}^{n \times n}$. 
The context-sensitive notation $|\bullet|$ may also denote the absolute value of a scalar, the Euclidean norm of a vector, or the Lebesgue measure of a set.
For any symmetric matrices $A , B \in \mathbb{S}$, $A \leq B$ means that all eigenvalues of $B - A \in \mathbb{S}$ are nonnegative.
The notation $A \lesssim B$ abbreviates $A \leq CB$ for a generic constant $C$ independent of the mesh-size and $A \approx B$ abbreviates $A \lesssim B \lesssim A$.
An open set $\omega \subset \mathbb{R}^n$ with boundary $\partial \omega$ satisfies a uniform exterior cone condition with the (closed) cone $K$ if, for all $x \in \partial \Omega$, there exists a cone $K_x$ with vertex $x$ such that $K_x$ is congruent to $K$ and $K_x \cap \overline{\Omega} = \{x\}$.

\section{Preliminary results from PDE theory}\label{sec:mathematical-setting}
Throughout this paper, we always assume that $L$ is a uniformly elliptic operator, i.e., there exist positive (ellipticity) constants $0 < \lambda \leq \Lambda$ such that the coefficient matrix $A$ satisfies
$\lambda \mathrm{I}_n \leq A \leq \Lambda \mathrm{I}_n$ a.e.~in $\Omega$.
The following maximum principle is a fundamental result in the analysis of strong solutions and  plays a major role in the design and analysis of the finite element schemes of this paper.
\begin{theorem}[ABP maximum principle]\label{thm:ABP-maximum-principle}
	Let $\Omega \subset \mathbb{R}^n$ be an open bounded set, $A \in L^\infty(\Omega;\mathbb{S})$ with $\lambda \mathrm{I}_n \leq A \leq \Lambda \mathrm{I}_n$, $b \in L^\infty(\Omega;\mathbb{R}^n)$, and $0 \leq c \in L^\infty(\Omega)$.
	There exists a constant $\cnst{cnst:ABP}$ depending on $n$, $\lambda$, $\|b\|_{L^\infty(\Omega)}$, and $\mathrm{diam}(\Omega)$ such that any strong solution $u \in C(\overline{\Omega}) \cap W^{2,n}_\mathrm{loc}(\Omega)$ to \eqref{def:PDE} satisfies \eqref{ineq:ABP-L-infty}.
\end{theorem}
\begin{proof}
	The proof can be found in \cite[Section 9.1]{GilbargTrudinger2001} with
	the constant
	\begin{align*}
		\cnst{cnst:ABP}^n \leq \mathrm{diam}(\Omega)^n\big(\mathrm{exp}\big(2^{n-2}\mathrm{diam}(\Omega)^n(1 + \|b\|_{L^\infty(\Omega)}^n/\mathcal{D})/(w_n n^n)\big)-1\big)/\mathcal{D}.
	\end{align*}
	Here, $w_n = \pi^{n/2}/\Gamma(n/2+1)$ with the gamma function $\Gamma$ and
	$\lambda^n \leq \mathcal{D} \leq \Lambda^n$ denotes the (essential) infimum of the determinant of $A$ over $\Omega$.
	In 2d, $w_2 = \pi$.
\end{proof}
Recall the norm $\|\bullet\|_V$ from \eqref{ineq:ABP-L-infty}. While the ABP maximum principle states that $\|u\|_{L^\infty(\Omega)} \leq \|u\|_{V}$,
we cannot expect the reverse bound $\|u\|_{V} \lesssim \|u\|_{L^\infty(\Omega)}$ in general. In fact, under additional assumptions, $\|u\|_{V}$ is an upper bound for the $H^2$ norm of $u$.
\begin{remark}[$H^2$ error control]
	Suppose that $n = 2$, $b = 0$, $c = 0$, and $g = 0$.
	If $\Omega$ is convex or the boundary $\partial \Omega$ of $\Omega$ is of class $C^{1,1}$, then the strong solution $u$ to \eqref{def:PDE} satisfies $u \in H^2(\Omega)$ with the estimate $\|u\|_{H^2(\Omega)} \lesssim \|f\|_{L^2(\Omega)} = \|u\|_V$ \cite{GilbargTrudinger2001,SmearsSueli2013}.
	On the other hand, a H\"older inequality leads to $\|u\|_V \leq \|A\|_{L^\infty(\Omega)}\|\D^2 u\|_{L^2(\Omega)}$ and so, we have the equivalence $\|u\|_{H^2(\Omega)} \approx \|u\|_V$ of norms.
\end{remark}
Note that, in general, the assumption $u \in C(\overline{\Omega}) \cap W^{2,n}_\mathrm{loc}(\Omega)$ in \Cref{thm:ABP-maximum-principle} cannot be relaxed by $u \in C(\overline{\Omega}) \cap W^{2,p}_\mathrm{loc}(\Omega)$ for some $p < n$ due to a result by Alexandrov \cite{Alexandrov1966}.
The ABP maximum principle leads to uniqueness of strong solutions, while existence can be rather involved.
If $A$ is continuous, then the existence of strong solutions can be established following \cite[Chapter 9]{GilbargTrudinger2001}.
(If the boundary $\partial \Omega$ is additionally of class $C^{1,1}$, then the global regularity $u \in W^{2,n}(\Omega)$ is guaranteed.)
Unfortunately, the situation is more complicated for merely bounded but possibly discontinuous coefficient $A \in L^\infty(\Omega;\mathbb{S})$.
In at least three space dimensions $n \geq 3$, the counterexamples from \cite{Nadirashvili1997,Safonov1999} show that a general existence and uniqueness theory for \eqref{def:PDE} cannot exist without additional assumptions on the coefficient $A$, e.g., if $A$ satisfies the Cordes condition and $p = 2$ \cite{Talenti1965}.
Therefore, the theory of this paper applies to the following case without a~priori information on the exact solution $u$.
\begin{theorem}[existence and uniqueness of strong solutions]\label{thm:existence}
	Suppose that the coefficients of $L$ satisfy \Cref{assumption:coefficients}.
	Given $f \in L^n(\Omega)$ and $g \in C(\partial\Omega)$, there exists a unique strong solution $u \in C(\overline{\Omega}) \cap W^{2,n}_\mathrm{loc}(\Omega)$ to \eqref{def:PDE}.
	For any open subset $\omega \Subset \Omega$, there exists a constant $\newcnst\label{cnst:W2n-interior}$ depending on $n$, $\lambda$, $\Lambda$, $\|b\|_{L^\infty(\Omega)}$, $\mathrm{diam}(\Omega)$, and $\mathrm{dist}(\omega, \partial \Omega)$ such that
	\begin{align}\label{ineq:W2n-estimate}
		\|u\|_{W^{2,n}(\omega)} \leq \cnst{cnst:W2n-interior}(\|u\|_{L^\infty(\partial \Omega)} + \|f\|_{L^n(\Omega)}).
	\end{align}
\end{theorem}
\begin{proof}
	It is known that any Lipschitz domain $\Omega$ satisfies an exterior cone condition \cite[Theorem 1.2.2.2]{Grisvard2011}.
	The existence of strong solutions is stated in \cite[Theorem 9.30]{GilbargTrudinger2001} even under the weaker assumption $A \in C(\Omega;\mathbb{S}) \cap L^\infty(\Omega;\mathbb{S})$ and the interior estimate \eqref{ineq:W2n-estimate} is given in \cite[Theorem 9.11]{GilbargTrudinger2001}. (Notice that the term $\|u\|_{L^p(\Omega)}$ therein can be replaced by $\|u\|_{L^\infty(\partial \Omega)}$ thanks to the ABP maximum principle from \Cref{thm:ABP-maximum-principle}.)
\end{proof}
An immediate consequence of theorems \ref{thm:ABP-maximum-principle} and \ref{thm:existence} is that all strong solutions to \eqref{def:PDE} form a Banach space. Recall $V$ from \eqref{def:V} and $\|\bullet\|_V$ from \eqref{ineq:ABP-L-infty}.
\begin{proposition}[$V$ is Banach space]
	Suppose that the coefficients of $L$ satisfy \Cref{assumption:coefficients}. Then $V$ is a Banach space endowed with the norm $\|\bullet\|_V$.
\end{proposition}
\begin{proof}
	We only prove completeness of $V$. Given any Cauchy sequence $(v_j)_{j \in \mathbb{N}_0}$ in $V$, the definition of $\|\bullet\|_V$ implies that $(v_j|_{\partial \Omega})_{j \in \mathbb{N}}$ resp.~$(L v_j)_{j \in \mathbb{N}}$ are Cauchy sequences in the Banach space $C(\partial \Omega)$ resp.~$L^n(\Omega)$. Therefore, there exist $f \in L^n(\Omega)$ and $g \in C(\partial \Omega)$ with $\lim_{j \to \infty} \|g - v_j\|_{L^\infty(\partial \Omega)} = 0$ and $\lim_{j \to \infty} \|f - Lv_j\|_{L^n(\Omega)} = 0$.
	\Cref{thm:existence} proves that there exists a unique strong solution $v \in C(\overline{\Omega}) \cap W^{2,n}_\mathrm{loc}(\Omega)$ to $L v = f$ in $\Omega$ and $v = g$ on $\partial \Omega$. In particular, $v \in V$ is the limit of the Cauchy sequence $(v_j)_{j \in \mathbb{N}}$ with respect to the norm $\|\bullet\|_V$.
\end{proof}
The next result states H\"older continuity of strong solutions to \eqref{def:PDE}.
\begin{theorem}[global H\"older regularity]\label{thm:Hoelder-regularity}
	Given $A \in L^\infty(\Omega;\mathbb{S})$ with $\lambda \mathrm{I}_n \leq A \leq \Lambda \mathrm{I}_n$, $b \in L^\infty(\Omega;\mathbb{R}^n)$, $0 \leq c \in L^\infty(\Omega)$, $f \in L^n(\Omega)$, and $g \in C^{0,\beta}(\partial \Omega)$ for some $\beta \in (0,1)$, then any strong solution $u \in C(\overline{\Omega}) \cap W^{2,n}_\mathrm{loc}(\Omega)$ to \eqref{def:PDE} is H\"older continuous $u \in C^{0,\alpha}(\overline{\Omega})$ with a positive parameter $\alpha \in (0,1)$ that solely depends on $n$, $\beta$, $\lambda$, $\Lambda$, $\|b\|_{L^\infty(\Omega)}$, and the cone condition of $\Omega$. In other words,
	\begin{align*}
		|u(x) - u(y)| \leq \cnst{cnst:Hoelder-regularity}|x - y|^\alpha
	\end{align*}
	for any $x,y \in \overline{\Omega}$. Here, the constant $\newcnst\label{cnst:Hoelder-regularity}$ solely depends on
	$n$, $\lambda$, $\Lambda$, $\|b\|_{L^\infty(\Omega)}$, $\|c\|_{L^\infty(\Omega)}$, $\|f\|_{L^n(\Omega)}$, $\|g\|_{C^{0,\beta}(\partial \Omega)}$, and the cone condition of $\Omega$.
\end{theorem}
\begin{proof}
	A proof of this result can be found in \cite[Theorem 6.2]{KoikeSwiech2009} or, in a slightly different formulation, in \cite[Chapter 9]{GilbargTrudinger2001}.
\end{proof}
In finite elements schemes, the coefficients of $L$ are approximated whenever numerical integration is used. Well-known results from \cite{Nadirashvili1997,Safonov1999} show that uniqueness may fail whenever we approximate general discontinuous coefficient $A$.
This issue does not arise if strong solutions in $C(\overline{\Omega}) \cap W^{2,n}_\mathrm{loc}(\Omega)$ exist.
\begin{lemma}[approximation of differential operator]\label{lem:approximation}
	Let $A \in L^\infty(\Omega;\mathbb{S})$ with $\lambda \mathrm{I}_n \leq A \leq \Lambda \mathrm{I}_n$, $b \in L^\infty(\Omega;\mathbb{R}^n)$, $0 \leq c \in L^\infty(\Omega)$, $f \in L^n(\Omega)$, $g \in C^{0,\beta}(\partial \Omega)$ and 
	$(A_j)_j \subset C(\Omega;\mathbb{S}) \cap L^\infty(\Omega;\mathbb{S})$, $(b_j)_j \subset L^\infty(\Omega;\mathbb{R}^n)$, $0 \leq (c_j)_j \subset L^\infty(\Omega)$, $(f_j)_j \subset L^n(\Omega)$, $(g_j)_j \subset C^{0,\beta}(\partial \Omega)$ for some $\beta \in (0,1)$ be given such that
	\begin{enumerate}[wide]
		\item[(a)] $A_j \to A$, $b_j \to b$, and $c_j \to c$ pointwise a.e.~in $\Omega$ as $j \to \infty$, $\|f - f_j\|_{L^n(\Omega)} \to 0$, and $\|g - g_j\|_{C^{0,\beta}(\partial \Omega)} \to 0$ as $j \to \infty$,
		\item[(b)] $\lambda \mathrm{I}_n \leq A_j \leq \Lambda \mathrm{I}_n$ pointwise in $\Omega$ for all $j \in \mathbb{N}$ and some constants $0 < \lambda \leq \Lambda$,
		\item[(c)] the $L^\infty$ norms of $b_j$ and $c_j$ are uniformly bounded independent of $j$.
	\end{enumerate}
	Suppose that there exists a strong solution $u \in C(\overline{\Omega}) \cap W^{2,n}_\mathrm{loc}(\Omega)$ to \eqref{def:PDE},
	then the sequence of strong solutions $u_j \in C(\overline{\Omega}) \cap W^{2,n}_\mathrm{loc}(\Omega)$ to
	\begin{align*}
		L_j u_j \coloneqq -A_j:\D^2 u_j + b_j \cdot \nabla u_j + c_j u_j = f_j \text{ in } \Omega \quad\text{and}\quad u_j = g_j \text{ on } \partial \Omega
	\end{align*}
	converges uniformly to $u$, i.e., $\lim_{j \to \infty} \|u - u_j\|_{L^\infty(\Omega)} = 0$.
\end{lemma}
We note that the existence of $u_j$ in \Cref{lem:approximation} follows from \Cref{thm:existence} because the leading coefficient $A_j$ is continuous.
For strong solutions $u \in W^{2,n}(\Omega)$, the assertion of \Cref{lem:approximation} can be found in \cite[Corollary 2.2]{Safonov1999}.
A proof under the assumption $u \in C(\overline{\Omega}) \cap W^{2,n}_\mathrm{loc}(\Omega)$ is given in \cite{CaffarelliCrandallKocanSwiech1996} even for fully nonlinear partial differential operators.
For the convenience of the reader, we provide an elementary proof for the linear case below.
\begin{proof}
	Let $v_j \coloneqq u - u_j$ be the strong solution to
	\begin{align*}
		L_j v_j = L_j(u - u_j) \text{ in } \Omega \quad\text{and}\quad v_j = g - g_j \text{ on } \partial \Omega.
	\end{align*}
	From \Cref{thm:Hoelder-regularity} and the assumptions (a)--(c), we deduce that $v_j$ is H\"older continuous with $\|v_j\|_{C^{0,\alpha}(\overline{\Omega})} \leq \cnst{cnst:proof-approximation}$ for some exponent $\alpha \in (0,1)$ and constant $\newcnst\label{cnst:proof-approximation}$ independent of the index $j$.
	Given $\varepsilon > 0$, define the open subset $\Omega_\varepsilon \coloneqq \{x \in \Omega: \mathrm{dist}(x,\partial \Omega) > (\varepsilon/\cnst{cnst:proof-approximation})^{1/\alpha}\} \Subset \Omega$.
	For any $x \in \Omega$ and $z \in \partial \Omega$, the H\"older regularity of $v_j$ proves		$|v_j(x)| \leq |v_j(x) - v_j(z)| + |v_j(z)| \leq C|x - z|^\alpha +  |g(z) - g_j(z)|$.
	This and the definition of $\Omega_\varepsilon$ imply
	\begin{align}\label{ineq:proof-approximation-error-near-boundary}
		\|v_j\|_{L^\infty(\Omega \setminus \Omega_\varepsilon)} = \|u - u_j\|_{L^\infty(\Omega\setminus\Omega_\varepsilon)} \leq \|g - g_j\|_{L^\infty(\partial \Omega)} + \varepsilon
	\end{align}
	The ABP maximum principle from \Cref{thm:ABP-maximum-principle} provides
	\begin{align}\label{ineq:proof-approximation-error-inside}
		\|v_j\|_{L^\infty(\Omega_\varepsilon)} = \|u-u_j\|_{L^\infty(\Omega_\varepsilon)} \leq \|v_j\|_{L^\infty(\partial \Omega_\varepsilon)} + \cnst{cnst:ABP}\|L_j(u - u_j)\|_{L^n(\Omega_\varepsilon)}.
	\end{align}
	A triangle inequality and $L u = f$ a.e.~in $\Omega$ lead to
	\begin{align*}
		&\|L_j(u - u_j)\|_{L^n(\Omega_\varepsilon)} \leq \|f - f_j\|_{L^n(\Omega)} + \|(L - L_j) u\|_{L^n(\Omega_\varepsilon)}.
	\end{align*}
	The combination of this with \eqref{ineq:proof-approximation-error-near-boundary}--\eqref{ineq:proof-approximation-error-inside} results in
	\begin{align*}
		&\|u - u_j\|_{L^\infty(\Omega)} \leq \|g - g_j\|_{L^\infty(\partial \Omega)} + \varepsilon + \cnst{cnst:ABP}(\|f - f_j\|_{L^n(\Omega)} + \|(L - L_j) u\|_{L^n(\Omega_\varepsilon)}).
	\end{align*}
	We observe that $(L - L_j) u \to 0$ pointwise a.e.~in $\Omega_\varepsilon$ as $j \to \infty$ from the assumption (a) and $(L - L_j) u \lesssim |\D^2 u| + |\nabla u| + |u|$ pointwise a.e.~in $\Omega_\varepsilon$ from the assumptions (b)--(c).
	Thus,	
	the Lebesgue dominated convergence theorem proves $\lim_{j \to \infty} \|(L - L_j) u\|_{L^n(\Omega_\varepsilon)} = 0$.
	Taking the limit of the previously displayed formula as $j \to \infty$ concludes $\limsup_{j \to \infty} \|u - u_j\|_{L^\infty(\Omega)} \leq \varepsilon$ for arbitrary $\varepsilon > 0$, whence $\lim_{j \to \infty} \|u - u_j\|_{L^\infty(\Omega)} = 0$.
\end{proof}
We note that the assumption $g \in C^{0,\beta}(\partial \Omega)$ in \Cref{lem:approximation} can be replaced by $g \in C(\partial \Omega)$ if $g_j = g$ for all $j$, i.e., if the Dirichlet data is not approximated.
The following density result is the foundation for the convergence analysis of this paper.
Recall the vector space $V$ from \eqref{def:V}.
\begin{lemma}[density]\label{thm:density}
	Suppose that the coefficients of $L$ satisfy \Cref{assumption:coefficients}.
	For any $v \in V$, there exists a sequence $(v_j)_j$ of functions $v_j \in W^{2,n}(\Omega)$ such that $L v_j = L v$ in $\Omega$ and $\lim_{j \to \infty} \|v - v_j\|_{L^\infty(\partial \Omega)} = 0$.
	In particular,
	$W^{2,n}(\Omega)$ is dense in $V$ (with respect to the norm $\|\bullet\|_V$).
\end{lemma}

\begin{proof}
	In the first step, the assertion is proven for any function $v \in V$ with homogenous boundary data $v = 0$ on $\partial \Omega$.
	Since $\Omega$ is Lipschitz, the set
	\begin{align*}
		\Omega(\delta) \coloneqq \{x \in \mathbb{R}^n : \mathrm{dist}(x,\Omega) < \delta\}
	\end{align*}
	is a Lipschitz domain for sufficiently small $0 < \delta \leq \delta_0$.
	In fact, the boundary of $\Omega(\delta)$ can be represented locally by the graph of some Lipschitz continuous function with the same Lipschitz constant in the same local coordinates as for $\Omega$ \cite[Theorem 4.1]{Doktor1976}.
	It is observed in \cite[p.~11]{Grisvard2011} that the cone condition of a Lipschitz domain solely depends on these parameters.
	Hence, the sequence $(\Omega_j)_j$ of Lipschitz domains
	\begin{align*}
		\Omega_j \coloneqq \{x \in \mathbb{R}^n : \mathrm{dist}(x,\Omega) < \delta_0/j\}
	\end{align*}
	approximates $\Omega$ with $\lim_{j \to \infty} \mathrm{dist}(\Omega,\partial \Omega_j) = 0$ and $\Omega_j$ satisfies a uniform exterior cone condition with a fixed cone $K$ independent of $j$.
	Let $A \in C(\mathbb{R}^n;\mathbb{S})$ be a (not relabelled) continuous extension of the coefficient $A$.
	In particular, $A$ is uniformly continuous in the compact set $\overline{\Omega}_1$.
	Therefore, there exists a $\gamma > 0$ such that $|A(x) - A(z)| \leq \lambda/2$ whenever $|x - z| \leq \gamma$ for all $x,z \in \overline{\Omega}_1$.
	The min-max principle shows, for any $x \in \Omega_1$ with $\mathrm{dist}(x,\partial \Omega) \leq \gamma$, that
	\begin{align}
		\min_{y \in S(\mathbb{R}^n)} y \cdot A(x) y \geq \min_{y \in S(\mathbb{R}^n)} y \cdot A(z) y - \min_{y \in S(\mathbb{R}^n)} y \cdot (A(z) - A(x)) y \geq \lambda/2,
	\end{align}
	where $z$ denotes the best-approximation of $x$ onto $\overline{\Omega}$ and $S(\mathbb{R}^n) = \{y \in \mathbb{R}^n : |y| = 1\}$.
	This shows $\lambda \mathrm{I}_n/2 \leq A$ and, by a similar argument, $A \leq (\Lambda + \lambda/2) \mathrm{I}_n$ in $\{z \in \overline{\Omega}_1: \mathrm{dist}(z,\Omega) \leq \gamma\}$. Without loss of generality we can assume that $\delta_0 \leq \gamma$ so that $\lambda \mathrm{I}_n/2 \leq A \leq (\Lambda + \lambda/2) \mathrm{I}_n$ holds pointwise in $\Omega_1$.
	(In particular, this holds in $\Omega_j$ for any $j \geq 1$.)
	For any $j \in \mathbb{N}$, let $v_j \in C(\overline{\Omega}_j) \cap W^{2,n}_\mathrm{loc}(\Omega_j)$ be the unique strong solution to
	\begin{align*}
		L v_j = f \text{ in } \Omega_j \quad\text{and}\quad v_j = 0 \text{ on } \partial \Omega_j,
	\end{align*}
	where the functions $b$, $c$, and $f$ are extended by zero outside $\Omega$.
	By design, $v_j|_{\Omega} \in W^{2,n}(\Omega)$ and we claim that $\lim_{j \to \infty} \|v - v_j\|_V = 0$. In fact,
	\begin{align}\label{eq:proof-density-V-norm}
		\|v - v_j\|_V = \|v - v_j\|_{L^\infty(\partial \Omega)} = \|v_j\|_{L^\infty(\partial \Omega)}.
	\end{align}
	From \Cref{thm:Hoelder-regularity}, we deduce that $v_j \in C^{0,\alpha}(\overline{\Omega}_j)$ with $\|v_j\|_{C^{0,\alpha}(\overline{\Omega}_j)} \leq \cnst{cnst:proof-density}$. The parameter $\alpha \in (0,1)$ and the constant $\newcnst\label{cnst:proof-density}$ are independent of $j$ because the cone condition of $\Omega_j$ is independent of $j$. The H\"older continuity of $v_j$ and $v_j = 0$ on $\partial \Omega_j$ provide
	\begin{align*}
		|v_j(x)| \leq \cnst{cnst:proof-density}\mathrm{dist}(x,\partial \Omega_j)^\alpha \quad\text{for any } x \in \partial \Omega.
	\end{align*}
	This and \eqref{eq:proof-density-V-norm} result in $\|v - v_j\|_{L^\infty(\partial \Omega)} \leq \cnst{cnst:proof-density} \mathrm{dist}(\Omega,\partial \Omega_j)^\alpha$, which tends to $0$ as $j \to \infty$.
	We thus proved that any $v \in V$ with $v = 0$ on $\partial \Omega$ can be approximated by functions in $W^{2,n}(\Omega)$.
	In the general case, let some $g_j \in C^\infty(\overline{\Omega})$ with $\|v - g_j\|_{L^\infty(\Omega)} \leq 1/(2j)$ for any $j \in \mathbb{N}$ be given. Then the strong solution $w \in C(\overline{\Omega}) \cap W^{2,n}_\mathrm{loc}(\Omega)$ to $L w = f - L g_j$ in $\Omega$ and $w = 0$ on $\partial \Omega$ satisfies $\|v - (w + g_j)\|_V = \|v - (w + g_j)\|_{L^\infty(\partial \Omega)} \leq 1/(2j)$.
	From the first step, there exists a $w_j \in W^{2,n}(\Omega)$ such that $L w_j = L w$ in $\Omega$ and $\|w - w_j\|_{L^\infty(\partial \Omega)} \leq 1/(2j)$. This and a triangle inequality conclude, for $v_j \coloneqq w_j + g_j \in W^{2,n}(\Omega)$, that 
	$L v_j = f$ in $\Omega$ and
	$\|v - v_j\|_V = \|v - v_j\|_{L^\infty(\partial \Omega)} \leq \|v - (w + g_j)\|_{L^\infty(\partial \Omega)} + \|w - w_j\|_{L^\infty(\partial \Omega)} \leq 1/j \to 0$ as $j \to \infty$.
\end{proof}
The following counterexample shows that the density result in \Cref{thm:density} may fail if we enforce Dirichlet boundary data pointwise onto the spaces therein.
\begin{proposition}[Laplace equation]\label{prop:failure-of-density}
	Let $u \in C(\overline{\Omega}) \cap H^2_\mathrm{loc}(\Omega)$ denote the strong solution to the Laplace problem $-\Delta u = 1$ in the two-dimensional $L$-shaped domain $\Omega \coloneqq (-1,1)^2 \setminus ([0,1] \times [-1,0])$ with homogenous boundary data $u = 0$ on $\partial \Omega$.
	Then $u$ cannot be the uniform limit of any bounded sequence $(u_j)_j \subset V$ of functions $u_j \in H^2(\Omega)$ with homogenous boundary data $u_j = 0$ on $\partial \Omega$.
	(Here, the boundedness of $(u_j)$ is understood with respect to the norm $\|\bullet\|_V$ from \eqref{ineq:ABP-L-infty}.)
\end{proposition}
\begin{proof}
	Since $f$ is smooth, $u \in C(\overline{\Omega}) \cap C^\infty(\Omega)$ \cite[Theorem 6.17]{GilbargTrudinger2001}. However, $u \notin H^2(\Omega)$ due to the reentrant corner of the domain $\Omega$.
	We recall the $H^2$ a~priori estimate $\|u_j\|_{H^2(\Omega)} \lesssim \|u_j\|_{L^2(\Omega)} + \|\Delta u_j\|_{L^2(\Omega)}$ for any $j \in \mathbb{N}$ on polygons from \cite[Theorem 4.3.1.4]{Grisvard2011}.
	This and the ABP maximum principle provide
	\begin{align*}
		\|u_j\|_{H^2(\Omega)} \lesssim \|\Delta u_j\|_{L^2(\Omega)}.
	\end{align*}
	In particular, $(u_j)$ is a bounded sequence with respect to the $H^2$ norm.
	The Banach-Alaoglu theorem proves that $u_j$ converges, up to some not relabelled subsequence, weakly to a $v \in H^2(\Omega)$.
	Since $u \notin H^2(\Omega)$, $v \neq u$.
	From the compact embedding $H^2(\Omega) \Subset C^{0,\alpha}(\overline{\Omega})$ for any $0 < \alpha < 1$ \cite[Theorem 6.3 III]{AdamsFournier2003}, we deduce that $u_j$ converges uniformly to $v$ up to another subsequence.
	Hence, the solution $u$ does not coincide with any accumulation point of $(u_j)_j$ with respect to the maximum norm. We note that this also holds for accumulation points with respect to the norm $\|\bullet\|_{W^{1,p}(\Omega)}$ for any $1 \leq p < \infty$ thanks to the compact embedding $H^2(\Omega) \Subset W^{1,p}(\Omega)$ \cite[Theorem 6.3 I]{AdamsFournier2003}.
\end{proof}

\section{Finite element approximation}\label{sec:convergencent-FEM}
Before the density result from \Cref{thm:density} is applied to the design of FEM, we fix some notation on the discrete level. Throughout the remaining parts of this paper, let $\Omega$ be a bounded polyhedral Lipschitz domain.

\subsection{Discrete spaces}
Let a quasi-uniform sequence $(\mathcal{T}_j)_j$ of regular triangulation of $\Omega$ into closed simplices or rectangles with the maximal mesh-size $h_j \coloneqq \max_{T \in \mathcal{T}_j} h_T$, where $h_T \coloneqq \mathrm{diam}(T)$ is the diameter of $T \in \mathcal{T}_j$, be given such that $\lim_{j \to \infty} h_j = 0$.
The set of all (resp.~interior and boundary) sides of $\mathcal{T}_j$ is denoted by $\mathcal{F}_j$ (resp.~$\mathcal{F}_j^i$ and $\mathcal{F}_j^b$).
For any interior side $F \in \mathcal{F}_j^i$, there exist two cells $T_+, T_- \in \mathcal{T}_j$ with $F = T_+ \cap T_-$.
The jump $[v]_F$ of any function $v \in W^{1,1}(T_\pm)$ is defined by $[v]_F \coloneqq v_{T_+} - v_{T_-}$.
Given any $T \in \mathcal{T}_j$ with sides $\mathcal{F}_j(T)$, $P_k(T)$ is the space of polynomials of degree at most $k \in \mathbb{N}$. The piecewise version of this reads $P_k(\mathcal{T}_j) \coloneqq \{v_j \in L^\infty(\Omega) : v_j|_T \in P_k(T) \text{ for all } T \in \mathcal{T}_j\}$.
Let $W^{\ell,p}(\mathcal{T}_j) \coloneqq \{v \in L^p(\Omega) : v|_T \in W^{\ell,p}(T) \text{ for all } T \in \mathcal{T}_j\}$, $\ell \geq 1$, $p \in [1,\infty]$, denote the space of piecewise $W^{\ell,p}$ functions, endowed with the norm
\begin{align*}
	\|v\|_{W^{\ell,p}(\mathcal{T}_j)} \coloneqq \Big(\sum\nolimits_{T \in \mathcal{T}_j} \|v\|_{W^{\ell,p}(T)}^p\Big)^{1/p},
\end{align*}
and $\nabla_\mathrm{pw} v$ (resp.~$\D^2_\mathrm{pw}$) denotes the piecewise gradient (resp.~Hessian) of $v \in W^{1,1}(\mathcal{T}_j)$ (resp.~$v \in W^{2,1}(\mathcal{T}_j)$) without explicit reference to the triangulation $\mathcal{T}_j$.

\subsection{Conforming FEM}\label{sec:conforming}
In this section, let $V(\mathcal{T}_j) \subset W^{2,\infty}(\Omega)$ be a $C^1$ conforming finite element space, e.g., the Argyris or Bogner-Fox-Schmit (BFS) finite element \cite{Ciarlet1978}. We assume that any $v \in W^{2,n}(\Omega)$ can be approximated by a sequence $(v_j)_j$ of discrete functions $v_j \in V(\mathcal{T}_j)$ such that $\lim_{j \to \infty} \|v - v_j\|_{W^{2,n}(\Omega)} = 0$.
The following result is an immediate consequence of \Cref{thm:density}.

\begin{corollary}[convergence of idealized FEM]\label{cor:convergence-FEM}
	Suppose that the coefficients of $L$ satisfy \Cref{assumption:coefficients}.
	Given $f \in L^n(\Omega)$ and $g \in C(\partial \Omega)$, any sequence $(u_j)_j$ of 
	discrete minimizers $u_j$ of the functional $\Phi$ from \eqref{def:Phi} in $V(\mathcal{T}_j)$
	converges uniformly to the strong solution $u \in C(\overline{\Omega}) \cap W^{2,n}_\mathrm{loc}(\Omega)$ to \eqref{def:PDE} as $j \to \infty$.
\end{corollary}
\begin{proof}
	Since $\|u - u_j\|_{L^\infty(\Omega)} \leq \|u - u_j\|_V = \Phi(u_j)$ from the ABP maximum principle in \Cref{thm:ABP-maximum-principle}, it suffices to show $\lim_{j \to \infty} \Phi(u_j) = 0$ for the convergence of FEM.
	Given $\varepsilon > 0$, \Cref{thm:density} proves that there exists a $v \in W^{2,n}(\Omega)$ such that $L v = f$ and $\|u - v\|_{L^\infty(\partial \Omega)} \leq \varepsilon/2$. Let $v_j$ be the best-approximation of $v$ in $V(\mathcal{T}_j)$ with respect to the $W^{2,n}$ norm.
	The triangle inequality provides
	\begin{align}\label{ineq:proof-conv-FEM-triangle-inequality}
		\Phi(u_j) \leq \Phi(v_j) \leq \|u - v\|_{V} + \|v - v_j\|_V \leq \varepsilon/2 + \|v - v_j\|_V.
	\end{align}
	Due to the Sobolev embedding \cite[Theorem 4.12 II]{AdamsFournier2003}, there exists a constant $\newcnst\label{cnst:Sobolev-embedding}$ depending on the domain $\Omega$ such that $\|w\|_{L^\infty(\Omega)} \leq \cnst{cnst:Sobolev-embedding}\|w\|_{W^{2,n}(\Omega)}$ for any $w \in W^{2,n}(\Omega)$.
	This, the H\"older, and a Cauchy inequality lead to
	\begin{align}
		\|v - v_j\|_V = \|v - v_j\|_{L^\infty(\Omega)} + \cnst{cnst:ABP}\|L(v - v_j)\|_{L^n(\Omega)}&\nonumber\\
		\leq \cnst{cnst:Sobolev-embedding}\|v - v_j\|_{W^{2,n}(\Omega)}
		+\,\cnst{cnst:ABP}\cnst{cnst:Hoelder-Cauchy}\|v - v_j\|_{W^{2,n}(\Omega)}&
		\label{ineq:proof-convergence-FEM-Sobolev-embedding}
	\end{align}
	with the constant $\newcnst\label{cnst:Hoelder-Cauchy} \coloneqq \big(\|A\|_{L^\infty(\Omega)}^{n/(n-1)} + \|b\|_{L^\infty(\Omega)}^{n/(n-1)} + \|c\|_{L^\infty(\Omega)}^{n/(n-1)}\big)^{(n-1)/n}$.
	Since $\lim_{j \to \infty} \|v - v_j\|_{W^{2,n}(\Omega)} = 0$, the index $j$ can be chosen sufficiently large so that $\|v - v_j\|_V \leq \varepsilon/2$. This and \eqref{ineq:proof-conv-FEM-triangle-inequality} result in $\Phi(u_j) \leq \varepsilon$ for sufficiently large $j$, which concludes the assertion.
\end{proof}
Notice that $u_j$ from \Cref{cor:convergence-FEM} is a best-approximation of $u$ in the discrete space $V(\mathcal{T}_j)$ with respect to the norm $\|\bullet\|_V$ (although the uniqueness of $u_j$ cannot be guaranteed).
However, the computation of $u_j$ involves a non-smooth nonlinear minimization problem.
We avoid this by enforcing the nonsmooth boundary residual as linear side constrains -- a well-known approach from mathematical optimization.
(Recall that, in general, the side constraints cannot be avoided by enforcing appropriate boundary data on the finite element functions as shown in \Cref{prop:failure-of-density}.)
For simplicity, suppose that there exists a set of points $\mathcal{L}_{j,b} \subset \partial \Omega$ on the boundary such that the following estimate holds
\begin{align}
	\|v_j\|_{L^\infty(F)} \leq \cnst{cnst:norm-equivalence-boundary}\max_{z \in \mathcal{L}_{j,b} \cap F} |v_j(z)|
	\label{ineq:assumtion-boundary-discrete-function}
\end{align}
for any discrete function $v_j \in V(\mathcal{T}_j)$ and $F \in \mathcal{F}_j^b$ with a constant $\newcnst\label{cnst:norm-equivalence-boundary}$ independent of the index $j$.
(For example, if $V(\mathcal{T}_j) \coloneqq Q_3(\mathcal{T}_j) \cap W^{2,\infty}(\Omega)$ is the space of BFS finite element functions in 2d \cite{Ciarlet1978} -- this is the space of all global $C^1(\overline{\Omega})$ functions that are bicubic when restricted to any rectangle $T \in \mathcal{T}_j$ -- then we can choose $\mathcal{L}_{j,b}$ as the set of all Lagrange points associated with $P_3(F)$ for some $F \in \mathcal{F}_j^b$.)
From \eqref{ineq:assumtion-boundary-discrete-function}, we deduce that
\begin{align}\label{ineq:equivalence-norm-boundary}
	\|g_j - v_j\|_{L^\infty(\partial \Omega)} \leq \cnst{cnst:norm-equivalence-boundary} \max_{z \in \mathcal{L}_{j,b}} |g_j(z) - v_j(z)|
\end{align}
Given an approximation $g_j \in V(\mathcal{T}_j)$ of the boundary data $g$, we define the set
\begin{align*}
	\mathcal{A}(g_j,\mathcal{T}_j) \coloneqq \{(v_j, t) \in V(\mathcal{T}_j) \times \mathbb{R}_{\geq 0} : - t \leq g_j(z) - v_j(z) \leq t \text{ for all } z \in \mathcal{L}_{j,b}\}
\end{align*}
of admissible discrete functions.
For a fixed constant $\alpha > 0$,
the proposed finite element scheme minimizes
\begin{align}\label{def:F}
	\Psi(v_j, t) \coloneqq \alpha t^n + \|f - L v_j\|_{L^n(\Omega)}^n \quad\text{among } (v_j,t) \in \mathcal{A}(g_j,\mathcal{T}_j).
\end{align}
\begin{theorem}[convergence of conforming FEM]\label{thm:convergence-conforming-FEM}
	Suppose that the coefficients of $L$ satisfy \Cref{assumption:coefficients}. 
	Given $f \in L^n(\Omega)$, and $g \in C(\partial \Omega)$, let $(g_j)_j$ with $g_j \in V(\mathcal{T}_j)$ approximate $g$ on the boundary with $\lim_{j \to \infty} \|g - g_j\|_{L^\infty(\partial \Omega)} = 0$.
	For any fixed $\alpha > 0$, the minimum of the functional $\Psi$ from \eqref{def:F} in $\mathcal{A}(g_j,\mathcal{T}_j)$ is attained at some $(u_j, t_j) \in \mathcal{A}_j(g_j,\mathcal{T}_j)$ and vanishes in the limit as $j \to \infty$. The sequence $(u_j)_j$ satisfies
	\begin{align}
		\lim_{j \to \infty} \|u - u_j\|_{L^\infty(\Omega)} = 0.
		\label{ineq:convergence-mFEM}
	\end{align}
\end{theorem}
\begin{proof}
	Fix $j \geq 0$.
	Since $\mathcal{A}(g_j,\mathcal{T}_j)$ is a subset of $V(\mathcal{T}_j) \times \mathbb{R}$ defined by linear side constraints, it is convex and closed but not bounded.
	Let $(v_{j\ell}, t_{j\ell})_{\ell}$ denote an infimizing sequence of \eqref{def:F}. By definition of $\mathcal{A}(g_j,\mathcal{T}_j)$ and \eqref{ineq:equivalence-norm-boundary}, $\|g_j - v_{j\ell}\|_{L^\infty(\partial \Omega)} \leq \cnst{cnst:norm-equivalence-boundary}t_{j\ell}$. An immediate consequence of this and the infimizing property of $(v_{j\ell}, t_\ell)_{\ell}$ is the boundedness of $v_{j\ell}$ in $V(\mathcal{T}_j)$ (with respect to the $\|\bullet\|_V$ norm), whence $(v_{j\ell}, t_{j\ell})_\ell$ is bounded.
	Since $V(\mathcal{T}_j) \times \mathbb{R}$ is finite dimensional, an accumulation point of $(v_{j\ell}, t_{j\ell})_\ell$ in $\mathcal{A}(g_j,\mathcal{T}_j)$ exists, which is a discrete of minimizer of \eqref{def:F} due to the continuity of $\Psi$ in $V(\mathcal{T}_j) \times \mathbb{R}$.
	Let $(u_j,t_j)$ minimize $\Psi$ in $\mathcal{A}(g_j, \mathcal{T}_j)$.
	Given any $v_j \in V(\mathcal{T}_j)$, define $r_j \coloneqq \max_{z \in \mathcal{L}_{j,b}} |g_j(z) - v_j(z)|$ and
	\begin{align}\label{def:alt-psi}
		\widetilde{\Psi}(v_j; g_j) \coloneqq \alpha\max_{z \in \mathcal{L}_{j,b}} |g_j(z) - v_j(z)|^n + \|f - L v_j\|^n_{L^n(\Omega)}.
	\end{align}
	It holds $\Psi(u_j,t_j) \leq \Psi(v_j,r_j) = \widetilde{\Psi}(v_j;g_j)$ because $(v_j,r_j) \in \mathcal{A}(g_j,\mathcal{T}_j)$. 
	Thus, $u_j$ minimizes $\widetilde{\Psi}$ among $v_j \in V(\mathcal{T}_j)$.
	Since $\widetilde{\Psi}(v_j; g_j) \lesssim (\|g - g_j\|_{L^\infty(\partial \Omega)} + \|g - v_j\|_{L^\infty(\partial \Omega)} + \|f - L v_j\|_{L^n(\Omega)})^n$ for any $v_j \in V(\mathcal{T}_j)$ from \eqref{ineq:equivalence-norm-boundary} and a Jensen inequality,
	\begin{align}\label{ineq:proof-convergence-practical}
		\widetilde{\Psi}(u_j; g_j)^{1/n} = \min_{v_j \in V(\mathcal{T}_j)} \widetilde{\Psi}(v_j;g_j)^{1/n} \lesssim \|g - g_j\|_{L^\infty(\partial \Omega)} + \min_{v_j \in V(\mathcal{T}_j)} \Phi(v_j)
	\end{align}
	with $\Phi$ from \eqref{def:Phi}. This and \Cref{cor:convergence-FEM} yield  $\lim_{j \to \infty} \widetilde{\Psi}(u_j; g_j) = 0$.
	On the other hand, \Cref{thm:ABP-maximum-principle}, \eqref{ineq:equivalence-norm-boundary}, and a triangle inequality prove $\|u - u_j\|_{L^\infty(\Omega)} \leq \|g - g_j\|_{L^\infty(\partial \Omega)} + \|g_j - u_j\|_{L^\infty(\partial \Omega)} + \cnst{cnst:ABP}\|f - Lu_j\|_{L^n(\Omega)} \lesssim \|g - g_j\|_{L^\infty(\partial \Omega)} + \widetilde{\Psi}(u_j;g_j)^{1/n}$.
	The limit of this as $j \to \infty$ concludes \eqref{ineq:convergence-mFEM}.
\end{proof}
Under additional smoothness assumptions on the exact solution $u$, we can obtain the following a~priori estimate.
\begin{corollary}[a~priori for conforming FEM]\label{cor:a-priori-conforming}
	In the setting of \Cref{thm:convergence-conforming-FEM}, suppose that $u \in W^{2,n}(\Omega)$ and $g = g_j$ for any $j$. Then
	\begin{align}\label{ineq:a-priori-conform}
		\Psi(u_j, t_j) = \widetilde{\Psi}(u_j; g) \lesssim \min_{v_j \in V(\mathcal{T}_j)} \|u - v_j\|_{W^{2,n}(\Omega)}^n.
	\end{align}
\end{corollary}
\begin{proof}
	The equality in \eqref{ineq:a-priori-conform} follows from $t_j = \max_{z \in \mathcal{L}_{j,b}} |g(z) - u_j(z)|$.
	(Suppose otherwise, then $\widetilde{t}_j \coloneqq \max_{z \in \mathcal{L}_{j,b}} |g(z) - u_j(z)| < t_j$ and $(u_j,\widetilde{t}_j) \in \mathcal{A}(g,\mathcal{T}_j)$, but $\Psi(u_j,\widetilde{t}_j) < \Psi(u_j,t_j)$. This is a contradiction to the minimizing property of $(u_j,t_j)$.)
	The remaining inequality in \eqref{ineq:a-priori-conform} follows from \eqref{ineq:proof-convergence-practical} with $g = g_j$ and 	\eqref{ineq:proof-convergence-FEM-Sobolev-embedding}.
\end{proof}

\subsection{Nonconforming FEM}
This subsection proposes a nonconforming FEM on simplicial meshes in two or three space dimensions $n = 2,3$. Given $k \geq 2$, let $V_{\mathrm{nc}}(\mathcal{T}_j) \coloneqq P_k(\mathcal{T}_j)$ denote the discrete ansatz space. 
As outlined in the proof of \Cref{thm:convergence-conforming-FEM},
any boundary residual arising from the ABP maximum principle in \Cref{thm:ABP-maximum-principle} will be enforced as side constraints.
Let $\mathcal{L}_{j}^k$ denote the set of all Lagrange points associated with the splines $P_k(\mathcal{T}_j) \cap W^{1,\infty}(\Omega)$ \cite[Proposition 7.12]{ErnGuermond2021} and $\mathcal{L}_{j,b}^k \coloneqq \mathcal{L}_j^k \cap \partial \Omega$.
Given $g_j \in V_\mathrm{nc}(\mathcal{T}_j)$, we define the set
\begin{align}
	\mathcal{A}_\mathrm{nc}(g_j,\mathcal{T}_j) \coloneqq \{(v_j,t) \in V_\mathrm{nc}(\mathcal{T}_j) \times \mathbb{R}_{\geq 0} : -t \leq (g_j - v_j)|_T(z) \leq t &\nonumber\\
	\text{for all } T \in \mathcal{T}_j \text{ and } z \in \mathcal{L}_{j,b}^k \cap T\}&
	\label{def:A-nc}
\end{align}
of admissible discrete functions. 
The equivalence of norms in finite dimensional spaces leads to
a piecewise version of \eqref{ineq:equivalence-norm-boundary},
\begin{align}\label{ineq:equivalence-norm-boundary-nonconforming}
	\|g_j - v_j\|_{L^\infty(\partial \Omega)} \leq \cnst{cnst:norm-equivalence-boundary-nc} \max_{T \in \mathcal{T}_j} \max_{z \in \mathcal{L}_{j,b}^k \cap T} (g_j|_T - v_j|_T)(z) \quad\text{for any } v_j \in V_\mathrm{nc}(\mathcal{T}_j)
\end{align}
with a positive constant $\newcnst\label{cnst:norm-equivalence-boundary-nc}$ that solely depends on the dimension $n$ and the polynomial degree $k$.
Given fixed positive parameters $\alpha, \sigma > 0$, the proposed nonconforming FEM minimizes the functional
\begin{align}\label{def:F-nc}
	\Psi_{\mathrm{nc}}(v_j, t) \coloneqq \alpha t^n + \|f - L_\mathrm{pw} v_j\|_{L^n(\Omega)}^n + \sigma\mathrm{s}_j(v_j),
\end{align}
among $(v_j,t) \in \mathcal{A}_{\mathrm{nc}}(g_j,\mathcal{T}_j)$,
where $L_\mathrm{pw} v_j \coloneqq -A:\D^2_\mathrm{pw} v_j + b \cdot \nabla_\mathrm{pw} v_j + c v_j$ is the piecewise application of the differential operator $L$ to $v_j$ and $\mathrm{s}_j(v_j) \coloneqq \sum_{T \in \mathcal{T}_j} \mathrm{s}_j(v_j;T)$ with
\begin{align}\label{def:stabilization}
	\mathrm{s}_j(v_j;T) \coloneqq \sum_{F \in \mathcal{F}_j^i \cap \mathcal{F}_j(T)} \big(h_F^{1-2n}\|[v_j]_F\|_{L^n(F)}^n +  h_F^{1-n} \|[\nabla_\mathrm{pw} v_j]_F\|_{L^n(F)}^n\big)
\end{align}
denotes the stabilization for all $v_j \in V_\mathrm{nc}(\mathcal{T}_j)$.
Since discontinuous functions do not satisfy the ABP maximum principle from \Cref{thm:ABP-maximum-principle}, we require a connection between the discrete space $V_\mathrm{nc}(\mathcal{T}_j)$ and $W^{2,n}(\Omega)$.
This is provided by a local averaging operator $\mathcal{J}_j : V_\mathrm{nc}(\mathcal{T}_j) \to P_m(\widehat{\mathcal{T}}_j) \cap W^{2,\infty}(\Omega)$ that maps $v_j \in V_\mathrm{nc}(\mathcal{T}_j)$ onto a $C^1$ conforming piecewise polynomial function $\mathcal{J}_j v_j$ of degree $m \geq k$ in a subtriangulation $\widehat{\mathcal{T}}_j$ of $\mathcal{T}_j$.
These spaces are known as the Hsieh–Clough–Tocher (HCT) macro element \cite{CloughTocher1965,WorseyFarin1987} and are available for arbitrary polynomial degree $m \geq 3$, cf.~\cite{DouglasDupontPercellScott1979} for 2d
and \cite{GuzmanLischkeNeilan2022} for 3d.
\begin{lemma}[enrichment operator]\label{lem:enrichment}
	Let $k \geq 2$ be given.
	There exists a linear operator $\mathcal{J}_j : V_\mathrm{nc}(\mathcal{T}_j) \to P_m(\widehat{\mathcal{T}}_j) \cap W^{2,\infty}(\Omega)$ for some $m \geq k$ such that, for all $v_j \in V_\mathrm{nc}(\mathcal{T}_j)$, $T \in \mathcal{T}_j$, and $p \in (1,\infty)$,
	\begin{align}
		h_T^{-2p}\|v_j - \mathcal{J}_j v_j\|_{L^p(T)}^p + h_T^{-p}\|\nabla (v_j - \mathcal{J}_j v_j)\|_{L^p(T)}^p + \|\D^2(v_j - \mathcal{J}_j v_j)\|_{L^p(T)}^p&\nonumber\\
		\leq \cnst{cnst:enrichment} \sum_{F \in \mathcal{F}_j^i, F \cap \partial T \neq \emptyset} \big(h_F^{1-2p}\|[v_j]_F\|^p_{L^p(F)} + h_F^{1-p}\|[\nabla_\mathrm{pw} v_j]_F\|^p_{L^p(F)}\big)&
		\label{ineq:local-enrichment-operator}
	\end{align}
	with a constant $\newcnst\label{cnst:enrichment}$ that solely depends on $n$, $p$, $k$, $m$, and the shape regularity of $\mathcal{T}_j$.
	Here, $\widehat{\mathcal{T}}_j$ denotes a subtriangulation of $\mathcal{T}_j$ such that $h_{\widehat{\mathcal{T}}_j} \approx h_{\mathcal{T}_j}$ a.e.~in $\Omega$ and the shape regularity of $\widehat{\mathcal{T}}_j$ depends exclusively on the shape regularity of $\mathcal{T}_j$.
\end{lemma}
\begin{proof}
	Local averaging techniques with the estimate \eqref{ineq:local-enrichment-operator} are well understood in the literature \cite{BrennerGudiSung2010,GeorgoulisHoustonVirtanen2011,Gallistl2015,CarstensenPuttkammer2023}; we refer to aforementioned articles for a precise definition of $\mathcal{J}_j$ with $m = 3$ and omit further details.
\end{proof}
We state the main result of this subsection.
\begin{theorem}[convergence of dG FEM]\label{thm:convergence-dG-FEM}
	Suppose that the coefficients of $L$ satisfy \Cref{assumption:coefficients} and $k \geq 2$. 
	Given $f \in L^n(\Omega)$ and $g \in C(\partial \Omega)$,
	let $(g_j)_j$ with $g_j \in V_\mathrm{nc}(\mathcal{T}_j)$ approximate $g$ on the boundary with $\lim_{j \to \infty} \|g - g_j\|_{L^\infty(\partial \Omega)} = 0$.
	For any fixed $\alpha, \sigma > 0$, the minimum of the functional $\Psi_\mathrm{nc}$ from \eqref{def:F-nc} in $\mathcal{A}_\mathrm{nc}(g_j,\mathcal{T}_j)$ is attained at some $(u_j,t_j) \in \mathcal{A}_\mathrm{nc}(g_j,\mathcal{T}_j)$ and vanishes in the limit as $j \to \infty$. The sequence $(u_j)_j$ satisfies
	\begin{align}
		\lim_{j \to \infty} \|u - u_j\|_{L^\infty(\Omega)} = 0.
		\label{ineq:convergence-mFEM-nc}
	\end{align}
\end{theorem}
\begin{proof}
	Fix $\alpha, \sigma > 0$.
	The existence of minimizers can follow the arguments in the conforming case and the fact that $\|u - v_j\|_{L^\infty(\Omega)} \lesssim \|g - g_j\|_{L^\infty(\partial \Omega)} + \Psi_{\mathrm{nc}}(v_j,t)$ for any $(v_j,t) \in \mathcal{A}_\mathrm{nc}(g_j,\mathcal{T}_j)$. (The proof of this is carried out in Step 3 below with the observation $\|g_j - v_j\|_{L^\infty(\Omega)} \leq \cnst{cnst:norm-equivalence-boundary-nc} t$).
	Let $(u_j,t_j)$ minimize $\Psi_{\mathrm{nc}}$ in $\mathcal{A}_\mathrm{nc}(g_j,\mathcal{T}_j)$.
	Given $\varepsilon > 0$, we select a $v \in W^{2,n}(\Omega)$ such that $L v = f$ and $\|u - v\|_{L^\infty(\partial \Omega)} \leq \varepsilon$ from \Cref{thm:density}.
	Let $v_j \coloneqq \Pi_{\mathcal{T}_j}^k v \in V_{\mathrm{nc}}(\mathcal{T}_j)$ denote the $L^2$ projection of $v$ onto $V_\mathrm{nc}(\mathcal{T}_j)$ and $r_j \coloneqq \max_{T \in \mathcal{T}_j} \max_{z \in \mathcal{L}_{j,b}^k \cap T} (g_j|_T - v_j|_T)(z)$. The remaining parts of the proof are divided into four steps.\medskip
	
	\emph{Step 1:} Prove $\limsup_{j \to \infty} r_j \leq \varepsilon$.
	The proof of this departs from the split
	\begin{align}\label{ineq:trace-discrete-function-split}
			\|g_j - v_j\|_{L^\infty(\partial \Omega)} &\leq \|g_j - g\|_{L^\infty(\partial \Omega)} + \|g - v\|_{L^\infty(\partial \Omega)} + \|v - v_j\|_{L^\infty(\partial \Omega)}.
		\end{align}
	We claim that $\lim_{j \to \infty} \|v - v_j\|_{L^\infty(\Omega)} = 0$. This can follow from density arguments outlined below for the sake of completeness. Given any $\delta > 0$,
	choose $w \in C^\infty(\overline{\Omega})$ such that $\|v - w\|_{L^\infty(\Omega)} \leq \delta$.
	The triangle inequality implies
	\begin{align*}
			\|v - v_j\|_{L^\infty(\Omega)} \leq \|v - w\|_{L^\infty(\Omega)} + \|(1 - \Pi_{\mathcal{T}_j}^k) w\|_{L^\infty(\Omega)} + \|\Pi_{\mathcal{T}_j}^k (w - v)\|_{L^\infty(\Omega)}.
		\end{align*}
	This, the approximation property $\|(1 - \Pi_{\mathcal{T}_j}^k) w\|_{L^\infty(\Omega)} \lesssim \|h_{\mathcal{T}_j}\nabla w\|_{L^\infty(\Omega)} \to 0$ as $j \to \infty$ and the $L^\infty$ stability $\|\Pi_{\mathcal{T}_j}^k (w - v)\|_{L^\infty(\Omega)} \lesssim \|v - w\|_{L^\infty(\Omega)}$ of the $L^2$ projection \cite[Lemma 11.18]{ErnGuermond2021} result in $\limsup_{j \to \infty} \|v - v_j\|_{L^\infty(\Omega)} \leq \newcnst\label{cnst:proof-convergence-dG-FEM}\delta$, where the constant \cnst{cnst:proof-convergence-dG-FEM} is independent of $\delta$. Since $\delta$ can be chosen arbitrary, this provides $\lim_{j \to \infty} \|v - v_j\|_{L^\infty(\Omega)} = 0$.
	In combination with $\lim_{j \to \infty} \|g - g_j\|_{L^\infty(\partial \Omega)} = 0$ and $\|g - v\|_{L^\infty(\partial \Omega)} = \|u - v\|_V \leq \varepsilon$, we deduce from \eqref{ineq:trace-discrete-function-split} that
	$$
		\limsup_{j \to \infty} r_j \leq \limsup_{j \to \infty} \|g_j - v_j\|_{L^\infty(\partial \Omega)} \leq \varepsilon.
	$$
	
	\emph{Step 2:} Prove $\lim_{j \to \infty} \Psi_\mathrm{nc}(u_j, t_j) = 0$.
	The H\"older and a Cauchy inequality show
	\begin{align}
		&\|f - L_\mathrm{pw} v_j\|_{L^n(\Omega)} = \|L_\mathrm{pw}(v - v_j)\|_{L^n(\Omega)}\nonumber\\
		&~\leq \big(\|A\|_{L^\infty(\Omega)}^{n/(n-1)} + \|b\|_{L^\infty(\Omega)}^{n/(n-1)} + \|c\|_{L^\infty(\Omega)}^{n/(n-1)}\big)^{(n-1)/n} \|(1 - \Pi_{\mathcal{T}_j}^k) v\|_{W^{2,n}(\mathcal{T}_j)}.\label{ineq:Hoelder-Cauchy}
	\end{align}
	Since $\lim_{j \to \infty} \|(1 - \Pi_{\mathcal{T}_j}^k) v\|_{W^{2,n}(\mathcal{T}_j)} = 0$, this implies
	\begin{align}\label{ineq:convergence-psi-first-term}
		\lim_{j \to \infty} \|f - L_\mathrm{pw} v_j\|_{L^n(\Omega)} = 0
	\end{align}
	and it remains to prove that $\lim_{j \to \infty} \mathrm{s}_j(v_j) = 0$ for the stabilization $\mathrm{s}_j$ from \eqref{def:stabilization} of $v_j$.
	For any interior side $F \in \mathcal{F}_j^i$ with the neighbouring cells $T_+, T_- \in \mathcal{T}_j$ and $F = T_+ \cap T_-$, $[v]_F = 0$ and $[\nabla v]_F = 0$ (in the sense of traces). A triangle and a trace inequality imply $\|[v_j]_F\|_{L^n(F)} \leq h_F^{-1/n} \|v - v_j\|_{L^n(\omega_F)} + h_F^{(n-1)/n}\|\nabla (v - v_j)\|_{L^n(\omega_F)}$ and $\|[\nabla_\mathrm{pw} v_j]_F\|_{L^n(F)} \leq h_F^{-1/n} \|\nabla_\mathrm{pw} (v - v_j)\|_{L^n(\omega_F)} + h_F^{(n-1)/n}\|\D_\mathrm{pw}^2 (v - v_j)\|_{L^n(\omega_F)}$ with $\omega_F \coloneqq \mathrm{int}(T_+ \cup T_-)$.
	This and the approximation property of the $L^2$ projection $\Pi_{\mathcal{T}_j}^k$
	\cite[Lemma 11.18]{ErnGuermond2021} verify
	\begin{align}
		\mathrm{s}_j(v_j) \lesssim \|h_{\mathcal{T}_j}^{-2}(1 - \Pi_{\mathcal{T}_j}^k) v\|^n_{L^n(\Omega)} + \|h_{\mathcal{T}_j}^{-1}\nabla_\mathrm{pw}(1 - \Pi_{\mathcal{T}_j}^k)v\|_{L^n(\Omega)}^n\nonumber&\\
		+\,\|\D^2_\mathrm{pw}(1 - \Pi_{\mathcal{T}_j}^k) v\|_{L^n(\Omega)}^n \lesssim \|\D^2_\mathrm{pw}(1 - \Pi_{\mathcal{T}_j}^k) v\|_{L^n(\Omega)}^n&.\label{ineq:stabilization}
	\end{align}
	Since $\lim_{j \to \infty} \|\D^2_\mathrm{pw}(1 - \Pi_{\mathcal{T}_j}^k) v\|_{L^n(\Omega)} = 0$, we deduce from \eqref{ineq:convergence-psi-first-term}, the definition of $\Psi_\mathrm{nc}$ in \eqref{def:F-nc}, and $\limsup_{j \to \infty} r_j \leq \varepsilon$ from Step 2 that $\limsup_{j \to \infty} \Psi_{\mathrm{nc}}(v_j,r_j) \leq \varepsilon^n$. This proves $\limsup_{j \to \infty} \Psi_{\mathrm{nc}}(u_j,t_j) \leq \limsup_{j \to \infty} \Psi_{\mathrm{nc}}(v_j,r_j) \leq \varepsilon^n$.
	Since $\varepsilon$ is chosen arbitrary, the claim follows.
	\medskip
	
%
	
	\emph{Step 3:} Prove $\|u - u_j\|_{L^\infty(\Omega)} \leq \|g - g_j\|_{L^\infty(\partial \Omega)} + \newcnst\label{cnst:a-posteriori-nc}\Psi_{\mathrm{nc}}(u_j,t_j)^{1/n}$ for some positive constant $\cnst{cnst:a-posteriori-nc}$ independent of $j$.
	Recall the local averaging operator $\mathcal{J}_j$ from \Cref{lem:enrichment}.
	The point of departure is the split
	\begin{align}\label{ineq:proof-convergence-dG-split}
		\|u - u_j\|_{L^\infty(\Omega)} \leq \|u - \mathcal{J}_j u_j\|_{L^\infty(\Omega)} + \|\mathcal{J}_j u_j - u_j\|_{L^\infty(\Omega)}.
	\end{align}
	The application of the ABP maximum principle from \Cref{thm:ABP-maximum-principle} to the difference $u - \mathcal{J}_j u_j \in C(\overline{\Omega}) \cap W^{2,n}_\mathrm{loc}(\Omega)$ and a triangle inequality lead to
	\begin{align}
		\|u - \mathcal{J}_j u_j\|_{L^\infty(\Omega)} \leq \|g - \mathcal{J}_j u_j\|_{L^\infty(\partial \Omega)} + \cnst{cnst:ABP}\|f - L \mathcal{J}_j u_j\|_{L^n(\Omega)}
		\leq \|g - u_j\|_{L^\infty(\partial \Omega)}&\nonumber\\
		+\,\|u_j - \mathcal{J}_j u_j\|_{L^\infty(\partial \Omega)} + \cnst{cnst:ABP}\|f - L_\mathrm{pw} u_j\|_{L^n(\Omega)} + \cnst{cnst:ABP}\|L_\mathrm{pw}(u_j - \mathcal{J}_j u_j)\|_{L^n(\Omega)}&.
	\end{align}
	The H\"older and a Cauchy inequality as in \eqref{ineq:proof-convergence-FEM-Sobolev-embedding} provide $\|L_\mathrm{pw}(u_j - \mathcal{J}_j u_j)\|_{L^n(\Omega)} \lesssim \|u_j - \mathcal{J}_j u_j\|_{W^{2,n}(\mathcal{T}_j)}$. This and \Cref{lem:enrichment} result in
	\begin{align}\label{ineq:bound-f-Luh}
		\|f - L_\mathrm{pw} u_j\|_{L^n(\Omega)} + \|L_\mathrm{pw}(u_j - \mathcal{J}_j u_j)\|_{L^n(\Omega)} \lesssim \Psi_{\mathrm{nc}}(u_j,t_j)^{1/n}.
	\end{align}
	The function $u_j - \mathcal{J}_j u_j$ is a piecewise polynomial in $\widehat{\mathcal{T}}_j$.
	Since the shape regularity of $\widehat{\mathcal{T}}_j$ only depends on the shape regularity of $\mathcal{T}$ (cf.~\cite{GuzmanLischkeNeilan2022,CarstensenPuttkammer2023} for the three dimensional case),	
	a scaling argument and \Cref{lem:enrichment} provide
	\begin{align}\label{ineq:proof-convergence-nc-J-max-error}
		\|u_j - \mathcal{J}_j u_j\|_{L^\infty(\Omega)} \lesssim h_j^{-1}\|u_j - \mathcal{J}_j u_j\|_{L^n(\Omega)} \lesssim h_j \mathrm{s}_j(u_j)^{1/n}.
	\end{align}
	The combination of this with \eqref{ineq:proof-convergence-dG-split}--\eqref{ineq:bound-f-Luh} results in
	\begin{align}\label{ineq:proof-convergence-nc-step-3}
		\|u - u_j\|_{L^\infty(\Omega)} \leq \|g - u_j\|_{L^\infty(\partial \Omega)} + \cnst{cnst:proof-nc-convergence}(\|f - L_\mathrm{pw} u_j\|_{L^n(\Omega)} + \mathrm{s}_j(u_j)^{1/n}).
	\end{align}
	with a positive constant $\newcnst\label{cnst:proof-nc-convergence}$ independent of $j$. 
	Observe that, similarly to the conforming case, $t_j = \max_{T \in \mathcal{T}_j} \max_{z \in \mathcal{L}_{j,b}^k \cap T} (g_j|_T - u_j|_T)(z)$. Hence, a triangle inequality and \eqref{ineq:equivalence-norm-boundary-nonconforming} provide
	$\|g - u_j\|_{L^\infty(\partial \Omega)} \leq \|g - g_j\|_{L^\infty(\partial \Omega)} + \cnst{cnst:norm-equivalence-boundary-nc}t_j$. This and \eqref{ineq:proof-convergence-nc-step-3} conclude the proof.
	\medskip
		
	\emph{Step 4:} Conclusion of the proof. Since $\lim_{j \to \infty} \|g - g_j\|_{L^\infty(\partial \Omega)} = 0$ by assumption and $\lim_{j \to \infty} \Psi_{\mathrm{nc}}(u_j,t_j) = 0$ from Step 2, the limit of the assertion in Step 3 as $j \to \infty$ concludes \eqref{ineq:convergence-mFEM-nc}.\qedhere
\end{proof}
As for the conforming FEM of \Cref{sec:conforming}, we obtain a~priori error estimates for $\Psi_{\mathrm{nc}}(u_j)$ under additional regularity assumptions on the exact solution $u$.
\begin{corollary}[a~priori for nonconforming FEM]\label{cor:a-priori-nc}
	In the setting of \Cref{thm:convergence-dG-FEM}, suppose that $u \in W^{2,n}(\Omega)$ and $g = g_j$ for any $j$. Then
	\begin{align}
		\Psi_\mathrm{nc}(u_j, t_j) \lesssim \|(1 - \Pi_{\mathcal{T}_j}^k) u\|_{W^{2,n}(\mathcal{T})}.
		\label{ineq:a-priori-nonconforming}
	\end{align}
\end{corollary}
\begin{proof}
	Let $v_j \coloneqq \Pi_{\mathcal{T}_j}^k u$.
	Since $g = g_j$, $r_j \coloneqq \max_{T \in \mathcal{T}_j} \max_{z \in \mathcal{L}_{j,b}^k \cap T} (g - v_j|_T)(z) \leq \|g - v_j\|_{L^\infty(\partial \Omega)}$
	and a triangle inequality prove
	\begin{align}\label{ineq:proof-a-priori-nc}
		\cnst{cnst:norm-equivalence-boundary-nc}^{-1} r_j \leq \|g - \mathcal{J}_j v_j\|_{L^\infty(\partial \Omega)} + \|\mathcal{J}_j v_j - v_j\|_{L^\infty(\partial \Omega)}.
	\end{align}
	The Sobolev embedding \cite[Theorem 6.3 III]{AdamsFournier2003} and a triangle inequality provide $\|g - \mathcal{J}_j v_j\|_{L^\infty(\partial \Omega)} \lesssim \|u - \mathcal{J}_j v_j\|_{W^{2,n}(\Omega)} \leq \|(1 - \Pi_{\mathcal{T}_j}^k) u\|_{W^{2,n}(\mathcal{T}_j)} + \|v_j - \mathcal{J}_j v_j\|_{W^{2,n}(\mathcal{T}_j)}$.
	Since $\|v_j - \mathcal{J}_j v_j\|_{W^{2,n}(\mathcal{T}_j)} \lesssim \mathrm{s}_j(v_j)^{1/n}$ from \Cref{lem:enrichment}, this, \eqref{ineq:proof-convergence-nc-J-max-error}, and \eqref{ineq:proof-a-priori-nc} show $r_j \lesssim \|(1 - \Pi_{\mathcal{T}_j}^k) u\|_{W^{2,n}(\mathcal{T}_j)} + \mathrm{s}_j(v_j)^{1/n}$.
	The assertion then follows from $\Psi_\mathrm{nc}(u_j,t_j) \leq \Psi_{\mathrm{nc}}(v_j,r_j)$, \eqref{ineq:Hoelder-Cauchy}, and \eqref{ineq:stabilization}.
\end{proof}
The following remark on another application of the density result from \Cref{thm:density} concludes this section.
\begin{remark}[least-squares]\label{rem:least-squares}
	Suppose that $\|\bullet\|$ is a norm in the Banach space $V$ from \eqref{def:V} such that $\|v\| \lesssim \|v\|_V$ for any $v \in V$, i.e., $\|\bullet\|$ is a weaker norm than $\|\bullet\|_V$ from \eqref{ineq:ABP-L-infty}.
	Recall the $C^1$ conforming finite element space $V(\mathcal{T}_j)$ from \Cref{sec:conforming}.
	We deduce from the proof of \Cref{cor:convergence-FEM} that any sequence $(u_j)_j$ of best-approximation $u_j$ of $u$ in $V(\mathcal{T}_j)$ with respect to the norm $\|\bullet\|$ satisfies $\lim_{j \to \infty} \|u - u_j\| = 0$.
	The choice $\|v\| \coloneqq \big(\|v\|_{L^2(\partial \Omega)}^2 + \|L v\|_{L^2(\Omega)}^2\big)^{1/2}$ is of particular interest because the best approximation $u_j$ of $u$ in $V(\mathcal{T}_j)$ with respect to $\|\bullet\|$ is the minimizer of the quadratic functional
	\begin{align}\label{def:Psi-LS}
		\Psi_\mathrm{LS}^2(v_j) \coloneqq \|g - v_j\|^2_{L^2(\partial \Omega)} + \|f - L v_j\|^2_{L^2(\Omega)} \text{ among } v_j \in V(\mathcal{T}_j).
	\end{align}
	Since $\Psi_\mathrm{LS}$ is strongly convex, the minimizer $u_j$ of $\Psi_\mathrm{LS}$ in $V(\mathcal{T}_j)$ is unique and satisfies the discrete Euler-Lagrange equations
	\begin{align}\label{eq:dELE-least-squares}
		\int_\Omega L u_j L v_j \,\mathrm{d}x + \int_{\partial \Omega} u_j v_j \,\mathrm{d}s = \int_\Omega f L v_j \,\mathrm{d}x + \int_{\partial \Omega} g v_j \,\mathrm{d}s.
	\end{align}
	While the previously proposed FEM need to solve a quadratic programming in 2d or a nonlinear convex minimization problem in 3d, this least-squares approach leads to a linear system of equations \eqref{eq:dELE-least-squares}.
	However, convergence can only established in the nonstandard norm $\|\bullet\|$ and control over the maximum norm is forfeited.
	Nevertheless, we can compute $\Phi(u_j)$ with $\Phi$ from \eqref{def:Phi} to check for uniform convergence a~posteriori.
	Thanks to the enrichment operator from \Cref{lem:enrichment}, we can extend the least-squares approach to nonconforming FEM as well.
	This leads to the minimization of the functional
	\begin{align*}
		\Psi_\mathrm{LS,nc}^2(v_j) \coloneqq \|g - v_j\|_{L^2(\partial \Omega)}^2 + \|f - L_\mathrm{pw} v_j\|^2_{L^2(\Omega)} + \widetilde{s}_j(v_j) \text{ among } v_j \in V_\mathrm{nc}(\mathcal{T}_j),
	\end{align*}
	where $\widetilde{s}_j$ is the quadratic version of $\mathrm{s}_j$ from \eqref{def:stabilization}.
	The minimizer $u_j$ of $\Psi_\mathrm{LS,nc}$ in $V_\mathrm{nc}(\mathcal{T}_j)$ is unique and convergence holds in the sense that $\lim_{j \to \infty} \Psi_\mathrm{LS,nc}(u_j) = 0$ as well as $\lim_{j \to \infty} \|u - \mathcal{J}_j u_j\| = 0$.
\end{remark}

\section{Numerical examples}\label{sec:numerical-examples}
This section presents results for four numerical benchmarks in two-dimensional domains with $b = 0$ and $c = 0$.

\subsection{Preliminary remarks}
We implement the following FEM.
The first method is the conforming BFS-FEM from \eqref{def:F} in \Cref{sec:conforming} with the BFS finite element space $V(\mathcal{T}_h) \coloneqq Q_3(\mathcal{T}_h) \cap W^{2,\infty}(\Omega)$ \cite{Ciarlet1978} on rectangular meshes $\mathcal{T}_h$ as ansatz space.
The second method is the nonconforming (NC) FEM from \eqref{def:F-nc} in \Cref{thm:convergence-dG-FEM} with the default parameter $\sigma = 1$.
Notice that the residuals $\Psi$ and $\Psi_{\mathrm{nc}}$ are not localizable because they involve terms related to the $L^\infty$ error on the boundary.
Thus, a straightforward adaptive mesh-refining strategy is not available.
In this paper, we suggest a penalization of the boundary residual by a weighted $L^2$ norm. To be precise, let $u_h$ denote the discrete solution to the corresponding finite element scheme. The adaptive computations utilize the refinement indicator
\begin{align*}
	\eta(T) \coloneqq \begin{cases}
		\|f - L u_h\|_{L^2(T)}^2 + \alpha\sum_{F \in \mathcal{F}_h(T), F \subset \partial \Omega} \|h_F^{\beta}(g - u_h)\|_{L^2(F)}^2 \text{ for the BFS-FEM},\\
		\|f - L_\mathrm{pw} u_h\|_{L^2(T)}^2 + \mathrm{s}_h(u_h;T) + \alpha\sum_{F \in \mathcal{F}_h(T), F \subset \partial \Omega} \|h_F^{\beta}(g - u_h)\|_{L^2(F)}^2
	\end{cases}
\end{align*}
for the NC-FEM
with a parameter $\beta$ that controls the magnitude of the penalization of the boundary residual. 
The default value of $\beta$ is set to $1$.
We utilize the D\"orfler marking strategy, i.e., at each refinement step, a subset $\mathcal{M} \subset \mathcal{T}$ with minimal cardinality is selected such that
\begin{align*}
	\sum\nolimits_{T \in \mathcal{T}} \eta(T) \leq \frac{1}{2}\sum\nolimits_{T \in \mathcal{M}} \eta(T).
\end{align*}
The convergence history plots display the quantities of interest against the number of degrees of freedom $\mathrm{ndof}$.
(Notice that $\mathrm{ndof} \approx h^{-2}_\mathrm{max}$ for uniform meshes.)
Solid lines in the convergence history plots indicate adaptive mesh-refinements, while dashed lines are associated with uniform mesh-refinements.
We recall from \Cref{thm:ABP-maximum-principle} that 
$$
\Phi(u_h) = \|g - u_h\|_{L^\infty(\partial \Omega)} + \cnst{cnst:ABP}\Psi(u_h) \geq \|u - u_h\|_{L^\infty(\Omega)}
$$ 
is a guaranteed upper bound (GUB) of the error $\|u - u_h\|_{L^\infty(\Omega)}$ for conforming and from \eqref{ineq:proof-convergence-nc-step-3} in the proof \Cref{thm:convergence-dG-FEM} that
$$
\Phi_{\mathrm{nc}}(u_h) \coloneqq (\|g - u_h\|_{L^\infty(\partial \Omega)}^2 + \|f - L_\mathrm{pw} u_h\|_{L^2(\Omega)}^2 + \mathrm{s}_h(u_h))^{1/2} \gtrsim \|u - u_h\|_{L^\infty(\Omega)}
$$
is an a~posteriori error estimate for nonconforming FEM.
Both error estimates hold for arbitrary discrete $u_h$, so it is applicable to inexact solve.

The quadratic optimization problems with linear side constrains proposed in this paper have been realized with the iterative solver \texttt{quadprog} from the MATLAB standard library. The parameters of \texttt{quadprog} are set to $\texttt{ConstraintTolerance} = \texttt{OptimalityTolerance} = \texttt{StepTolerance} = 10^{-14}$.

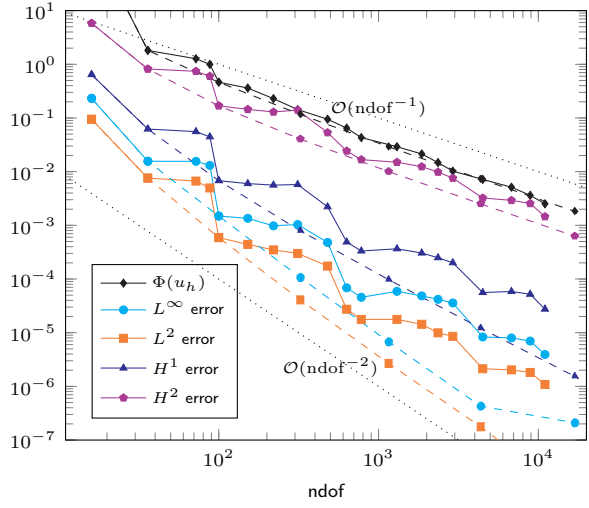
\begin{figure}
	\centering
	\begin{tikzpicture}
		\begin{loglogaxis}[legend pos=south west,legend cell align=left,
			legend style={fill=none},
			ymin=1e-7,ymax=1e1,
			xmin=11,xmax=2e4,
			ytick={1e-7,1e-6,1e-5,1e-4,1e-3,1e-2,1e-1,1e0,1e1,1e2},
			xtick={1e1,1e2,1e3,1e4,1e5,1e6,1e7}]
			\pgfplotsset{
				cycle list={%
					{Black, mark=diamond*, mark size=1.5pt},
					{Cyan, mark=*, mark size=1.5pt},
					{Orange, mark=square*, mark size=1.5pt},
					{Blue, mark=triangle*, mark size=1.5pt},
					{Mulberry, mark=pentagon*, mark size=1.5pt},
					{Black, dashed, mark=diamond*, mark size=1.5pt},
					{Cyan, dashed, mark=*, mark size=1.5pt},
					{Orange, dashed, mark=square*, mark size=1.5pt},
					{Blue, dashed, mark=triangle*, mark size=1.5pt},
					{Mulberry, dashed, mark=pentagon*, mark size=1.5pt}
				},
				legend style={
					at={(0.05,0.05)}, anchor=south west},
				font=\sffamily\scriptsize,
				xlabel=ndof,ylabel=
			}
			\addplot+ table[x=ndof,y=GUB]{convhist_ex_smooth_bfs_adap_alpha_1000_tau_2.dat};
			\addlegendentry{\tiny $\Phi(u_h)$}
			\addplot+ table[x=ndof,y=max_error]{convhist_ex_smooth_bfs_adap_alpha_1000_tau_2.dat};
			\addlegendentry{\tiny $L^\infty$ error}
			\addplot+ table[x=ndof,y=L2_error]{convhist_ex_smooth_bfs_adap_alpha_1000_tau_2.dat};
			\addlegendentry{\tiny $L^2$ error}
			\addplot+ table[x=ndof,y=H1_error]{convhist_ex_smooth_bfs_adap_alpha_1000_tau_2.dat};
			\addlegendentry{\tiny $H^1$ error}
			\addplot+ table[x=ndof,y=H2_error]{convhist_ex_smooth_bfs_adap_alpha_1000_tau_2.dat};
			\addlegendentry{\tiny $H^2$ error}
			\addplot+ table[x=ndof,y=GUB]{convhist_ex_smooth_bfs_unif_alpha_1000_tau_2.dat};
			\addplot+ table[x=ndof,y=max_error]{convhist_ex_smooth_bfs_unif_alpha_1000_tau_2.dat};
			\addplot+ table[x=ndof,y=L2_error]{convhist_ex_smooth_bfs_unif_alpha_1000_tau_2.dat};
			\addplot+ table[x=ndof,y=H1_error]{convhist_ex_smooth_bfs_unif_alpha_1000_tau_2.dat};
			\addplot+ table[x=ndof,y=H2_error]{convhist_ex_smooth_bfs_unif_alpha_1000_tau_2.dat};
			\addplot+ [dotted,mark=none,black] coordinates{(1e1,1e1) (1e5,1e-3)};
			\node(z) at  (axis cs:1e3,7e-2)
			[above] {\tiny $\mathcal O (\mathrm{ndof}^{-1})$};
			\addplot+ [dotted,mark=none,black] coordinates{(1e1,1e-2) (1e5,1e-10)};
			\node(z) at  (axis cs:5e2,1e-6)
			[above] {\tiny $\mathcal O (\mathrm{ndof}^{-2})$};
		\end{loglogaxis}
	\end{tikzpicture}
	\caption{Convergence history of the BFS-FEM for the first experiment with $\alpha = 10^3$ (solid lines = adaptive, dashed lines = uniform). \label{f:exp-smooth-conv-bfs}
	}
\end{figure}

\begin{figure}
	\centering
	\begin{subfigure}[b]{0.48\textwidth}
		\begin{tikzpicture}[scale=0.78]
			\begin{loglogaxis}[legend pos=south west,legend cell align=left,
				legend style={fill=none},
				ymin=1e-5,ymax=3e2,
				xmin=20,xmax=2e5,
				ytick={1e-7,1e-6,1e-5,1e-4,1e-3,1e-2,1e-1,1e0,1e1,1e2,1e3},
				xtick={1e1,1e2,1e3,1e4,1e5,1e6,1e7}]
				\pgfplotsset{
					cycle list={%
						{Black, mark=diamond*, mark size=1.5pt},
						{Cyan, mark=*, mark size=1.5pt},
						{Orange, mark=square*, mark size=1.5pt},
						{Blue, mark=triangle*, mark size=1.5pt},
						{Black, dashed, mark=diamond*, mark size=1.5pt},
						{Cyan, dashed, mark=*, mark size=1.5pt},
						{Orange, dashed, mark=square*, mark size=1.5pt},
						{Blue, dashed, mark=triangle*, mark size=1.5pt},
						{Mulberry, dashed, mark=pentagon*, mark size=1.5pt}
					},
					legend style={
						at={(0.05,0.05)}, anchor=south west},
					font=\sffamily\scriptsize,
					xlabel=ndof,ylabel=
				}
				\addplot+ table[x=ndof,y=eta]{convhist_ex_smooth_dG_adap_k_0_alpha_1000_tau_2.dat};
				\addlegendentry{\tiny $k = 0$}
				\addplot+ table[x=ndof,y=eta]{convhist_ex_smooth_dG_adap_k_1_alpha_1000_tau_2.dat};
				\addlegendentry{\tiny $k = 1$}
				\addplot+ table[x=ndof,y=eta]{convhist_ex_smooth_dG_adap_k_2_alpha_1000_tau_2.dat};
				\addlegendentry{\tiny $k = 2$}
				\addplot+ table[x=ndof,y=eta]{convhist_ex_smooth_dG_adap_k_3_alpha_1000_tau_2.dat};
				\addlegendentry{\tiny $k = 3$}
				\addplot+ table[x=ndof,y=eta]{convhist_ex_smooth_dG_unif_k_0_alpha_1000_tau_2.dat};
				\addplot+ table[x=ndof,y=eta]{convhist_ex_smooth_dG_unif_k_1_alpha_1000_tau_2.dat};
				\addplot+ table[x=ndof,y=eta]{convhist_ex_smooth_dG_unif_k_2_alpha_1000_tau_2.dat};
				\addplot+ table[x=ndof,y=eta]{convhist_ex_smooth_dG_unif_k_3_alpha_1000_tau_2.dat};
				\addplot+ [dotted,mark=none,black] coordinates{(1e1,4e2) (1e7,4e-1)};
				\node(z) at  (axis cs:3e3,1e1)
				[above] {\tiny $\mathcal O (\mathrm{ndof}^{-1/2})$};
				\addplot+ [dotted,mark=none,black] coordinates{(1e1,5e1) (1e5,5e-7)};
				\node(z) at  (axis cs:1.5e3,1e-3)
				[above] {\tiny $\mathcal O (\mathrm{ndof}^{-2})$};
			\end{loglogaxis}
		\end{tikzpicture}
	\end{subfigure}
	\begin{subfigure}[b]{0.48\textwidth}
		\begin{tikzpicture}[scale=0.78]
			\begin{loglogaxis}[legend pos=south west,legend cell align=left,
				legend style={fill=none},
				ymin=1e-7,ymax=1e1,
				xmin=20,xmax=3e5,
				ytick={1e-7,1e-6,1e-5,1e-4,1e-3,1e-2,1e-1,1e0,1e1,1e2},
				xtick={1e1,1e2,1e3,1e4,1e5,1e6,1e7}]
				\pgfplotsset{
					cycle list={%
						{Black, mark=diamond*, mark size=1.5pt},
						{Cyan, mark=*, mark size=1.5pt},
						{Orange, mark=square*, mark size=1.5pt},
						{Blue, mark=triangle*, mark size=1.5pt},
						{Black, dashed, mark=diamond*, mark size=1.5pt},
						{Cyan, dashed, mark=*, mark size=1.5pt},
						{Orange, dashed, mark=square*, mark size=1.5pt},
						{Blue, dashed, mark=triangle*, mark size=1.5pt},
						{Mulberry, dashed, mark=pentagon*, mark size=1.5pt}
					},
					legend style={
						at={(0.05,0.05)}, anchor=south west},
					font=\sffamily\scriptsize,
					xlabel=ndof,ylabel=
				}
				\addplot+ table[x=ndof,y=max_error]{convhist_ex_smooth_dG_adap_k_0_alpha_1000_tau_2.dat};
				\addlegendentry{\tiny $k = 0$}
				\addplot+ table[x=ndof,y=max_error]{convhist_ex_smooth_dG_adap_k_1_alpha_1000_tau_2.dat};
				\addlegendentry{\tiny $k = 1$}
				\addplot+ table[x=ndof,y=max_error]{convhist_ex_smooth_dG_adap_k_2_alpha_1000_tau_2.dat};
				\addlegendentry{\tiny $k = 2$}
				\addplot+ table[x=ndof,y=max_error]{convhist_ex_smooth_dG_adap_k_3_alpha_1000_tau_2.dat};
				\addlegendentry{\tiny $k = 3$}
				\addplot+ table[x=ndof,y=max_error]{convhist_ex_smooth_dG_unif_k_0_alpha_1000_tau_2.dat};
				\addplot+ table[x=ndof,y=max_error]{convhist_ex_smooth_dG_unif_k_1_alpha_1000_tau_2.dat};
				\addplot+ table[x=ndof,y=max_error]{convhist_ex_smooth_dG_unif_k_2_alpha_1000_tau_2.dat};
				\addplot+ table[x=ndof,y=max_error]{convhist_ex_smooth_dG_unif_k_3_alpha_1000_tau_2.dat};
				\addplot+ [dotted,mark=none,black] coordinates{(1e1,5e1) (1e6,5e-4)};
				\node(z) at  (axis cs:3e3,1e-1)
				[above] {\tiny $\mathcal O (\mathrm{ndof}^{-1})$};
				\addplot+ [dotted,mark=none,black] coordinates{(1e1,1e1) (1e5,1e-11)};
				\node(z) at  (axis cs:2e3,1e-6)
				[above] {\tiny $\mathcal O (\mathrm{ndof}^{-3})$};
			\end{loglogaxis}
		\end{tikzpicture}
	\end{subfigure}
	\caption{Convergence history of $\Phi_{\mathrm{nc}}(u_h)$ (left) and $L^\infty$ error (right) for the NC-FEM in the first experiment with various $k$ and $\alpha = 10^{3}$ (solid lines = adaptive, dashed lines = uniform).
		\label{f:exp-smooth-conv-dG}
	}
\end{figure}
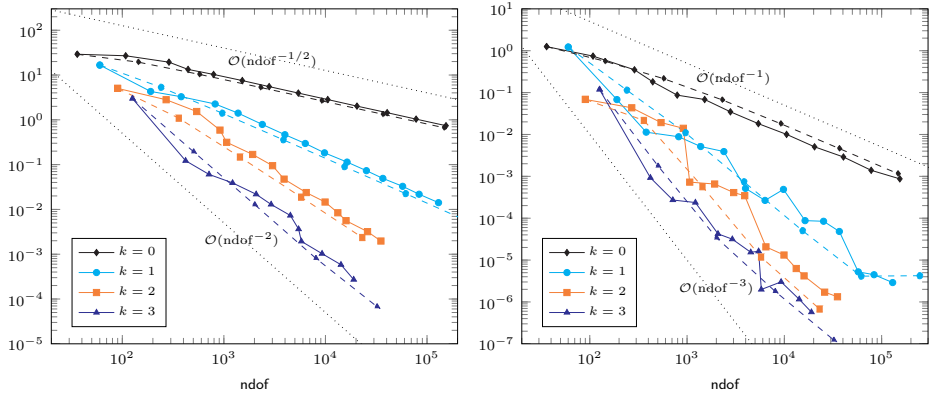

\subsection{First experiment}
This benchmark approximates the smooth exact solution
\begin{align*}
	u(x) \coloneqq \sin(\pi x_1)\sin(\pi x_2) \quad\text{for all } x = (x_1, x_2)
\end{align*}
to \eqref{def:PDE}
in the convex unit square $\Omega \coloneqq (0,1)^2$ with the coefficient matrix
$$A(r, \varphi) \coloneqq 
\begin{pmatrix}
	1 + r^{1/2} & - r^{1/2}\\ - r^{1/2} & 1 + 5r^{1/2}
\end{pmatrix}$$
in polar coordinates and the right-hand side $f(x) \coloneqq - A(x) : \D^2 u(x)$. Since $\Omega$ is convex and $A$ satisfies the Cordes condition, the framework of \cite{SmearsSueli2013} applies and leads to more efficient numerical methods. (This can be recovered by enforcing the Dirichlet boundary data strongly to the discrete ansatz space. It is straight-forward to verify that the resulting minimal residual method, which coincides with \cite{g2017b}, converges.)
While this application does not belong to the focus of this paper, it allows for an investigation of possible convergence rates in the smooth case.
Due to different scaling of norms in the objective functional \eqref{def:F} and \eqref{def:F-nc}, it is preferable to set $\alpha$ sufficiently large to counter numerical instabilities on fine meshes. In this smooth example, the errors are expected to become small and so, we set $\alpha \coloneqq 10^3$. \Cref{f:exp-smooth-conv-bfs}(a) displays the optimal convergence rate $1$ for $\Phi(u_h)$. This coincides with the a~priori result in \Cref{cor:a-priori-conforming}. The same convergence rate is observed for the $H^2$ error. The errors in the $L^\infty$, $L^2$, and $H^1$ norms converge with a faster convergence rate up to $2$ for the $L^2$ and $L^\infty$ errors. This provides empirical evidence that $\Phi(u_h)$ is not an efficient error estimator for the $L^\infty$ error in general.
Adaptive computations do not provide any improvements in this case.
\Cref{f:exp-smooth-conv-dG} displays the convergence history plot of the a~posteriori error estimator $\Psi_{\mathrm{nc}}(u_h)$ and the $L^\infty$ error for the NC-FEM with similar behavior to the conforming case: $\Psi_{\mathrm{nc}}(u_h)$ converges optimally with the convergence rates $(k+1)/2$ as predicted in \Cref{cor:a-priori-nc} and the $L^\infty$ error converges faster than $\Psi_{\mathrm{nc}}(u_h)$.
Undisplayed numerical results for different values of $\alpha$ show no changes in convergence rates of the displayed quantities.

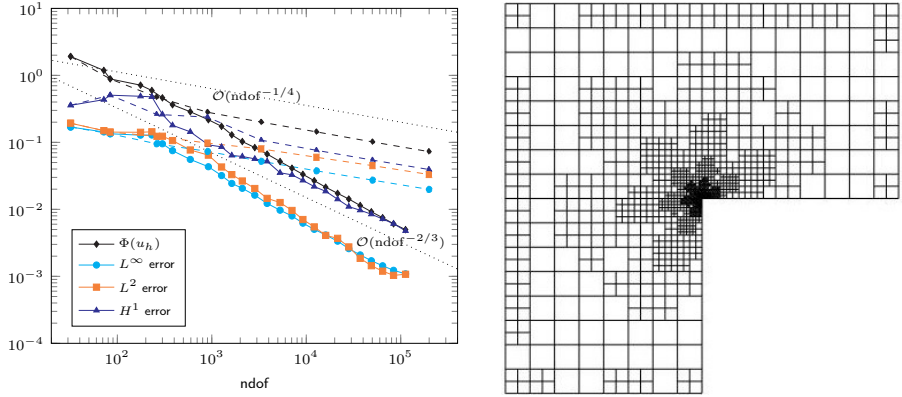
\begin{figure}
	\centering
	\begin{subfigure}[b]{0.48\textwidth}
		\begin{tikzpicture}[scale=0.78]
			\begin{loglogaxis}[legend pos=south west,legend cell align=left,
				legend style={fill=none},
				ymin=1e-4,ymax=1e1,
				xmin=20,xmax=4e5,
				ytick={1e-7,1e-6,1e-5,1e-4,1e-3,1e-2,1e-1,1e0,1e1,1e2},
				xtick={1e1,1e2,1e3,1e4,1e5,1e6,1e7}]
				\pgfplotsset{
					cycle list={%
						{Black, mark=diamond*, mark size=1.5pt},
						{Cyan, mark=*, mark size=1.5pt},
						{Orange, mark=square*, mark size=1.5pt},
						{Blue, mark=triangle*, mark size=1.5pt},
						{Black, dashed, mark=diamond*, mark size=1.5pt},
						{Cyan, dashed, mark=*, mark size=1.5pt},
						{Orange, dashed, mark=square*, mark size=1.5pt},
						{Blue, dashed, mark=triangle*, mark size=1.5pt},
						{Mulberry, dashed, mark=pentagon*, mark size=1.5pt}
					},
					legend style={
						at={(0.05,0.05)}, anchor=south west},
					font=\sffamily\scriptsize,
					xlabel=ndof,ylabel=
				}
				\addplot+ table[x=ndof,y=GUB]{convhist_ex1_bfs_adap_alpha_10_tau_2.dat};
				\addlegendentry{\tiny $\Phi(u_h)$}
				\addplot+ 	table[x=ndof,y=max_error]{convhist_ex1_bfs_adap_alpha_10_tau_2.dat};
				\addlegendentry{\tiny $L^\infty$ error}
				\addplot+ 	table[x=ndof,y=L2_error]{convhist_ex1_bfs_adap_alpha_10_tau_2.dat};
				\addlegendentry{\tiny $L^2$ error}
				\addplot+ 	table[x=ndof,y=H1_error]{convhist_ex1_bfs_adap_alpha_10_tau_2.dat};
				\addlegendentry{\tiny $H^1$ error}
				\addplot+ table[x=ndof,y=GUB]{convhist_ex1_bfs_unif_alpha_10_tau_2.dat};
				\addplot+ 	table[x=ndof,y=max_error]{convhist_ex1_bfs_unif_alpha_10_tau_2.dat};
				\addplot+ 	table[x=ndof,y=L2_error]{convhist_ex1_bfs_unif_alpha_10_tau_2.dat};
				\addplot+ 	table[x=ndof,y=H1_error]{convhist_ex1_bfs_unif_alpha_10_tau_2.dat};
				\addplot+ [dotted,mark=none,black] coordinates{(1e1,2e0) (1e9,2e-2)};
				\node(z) at  (axis cs:3e3,3e-1)
				[above] {\tiny $\mathcal O (\mathrm{ndof}^{-1/4})$};
				\addplot+ [dotted,mark=none,black] coordinates{(1e1,1.5e0) (1e10,1.5e-6)};
				\node(z) at  (axis cs:1e5,2e-3)
				[above] {\tiny $\mathcal O (\mathrm{ndof}^{-2/3})$};
			\end{loglogaxis}
		\end{tikzpicture}
	\end{subfigure}
	\begin{subfigure}[b]{0.48\textwidth}
		\centering
		\includegraphics[scale=0.57]{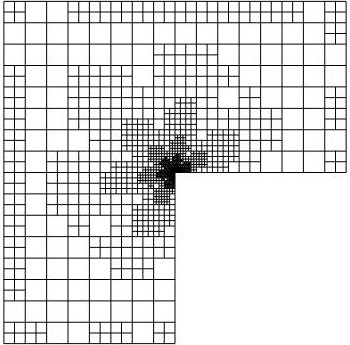}
	\end{subfigure}
	\caption{Convergence history (left) and adaptive mesh with 2619 elements (right) of the BFS-FEM for the second experiment with $\alpha = 10$ (solid lines = adaptive, dashed lines = uniform). \label{f:exp1-conv-bfs}
	}
\end{figure}

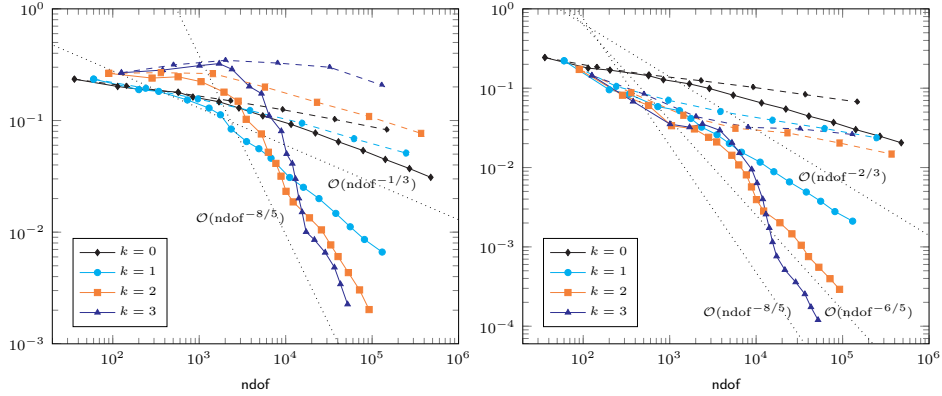
\begin{figure}
	\centering
	\begin{subfigure}[b]{0.48\textwidth}
		\begin{tikzpicture}[scale=0.78]
			\begin{loglogaxis}[legend pos=south west,legend cell align=left,
				legend style={fill=none},
				ymin=1e-3,ymax=1e0,
				xmin=20,xmax=1e6,
				ytick={1e-7,1e-6,1e-5,1e-4,1e-3,1e-2,1e-1,1e0,1e1,1e2,1e3},
				xtick={1e1,1e2,1e3,1e4,1e5,1e6,1e7}]
				\pgfplotsset{
					cycle list={%
						{Black, mark=diamond*, mark size=1.5pt},
						{Cyan, mark=*, mark size=1.5pt},
						{Orange, mark=square*, mark size=1.5pt},
						{Blue, mark=triangle*, mark size=1.5pt},
						{Black, dashed, mark=diamond*, mark size=1.5pt},
						{Cyan, dashed, mark=*, mark size=1.5pt},
						{Orange, dashed, mark=square*, mark size=1.5pt},
						{Blue, dashed, mark=triangle*, mark size=1.5pt},
						{Mulberry, dashed, mark=pentagon*, mark size=1.5pt}
					},
					legend style={
						at={(0.05,0.05)}, anchor=south west},
					font=\sffamily\scriptsize,
					xlabel=ndof,ylabel=
				}
				\addplot+ table[x=ndof,y=eta]{convhist_ex1_dG_adap_k_0_alpha_1_tau_2.dat};
				\addlegendentry{\tiny $k = 0$}
				\addplot+ table[x=ndof,y=eta]{convhist_ex1_dG_adap_k_1_alpha_10_tau_2.dat};
				\addlegendentry{\tiny $k = 1$}
				\addplot+ table[x=ndof,y=eta]{convhist_ex1_dG_adap_k_2_alpha_100_tau_2.dat};
				\addlegendentry{\tiny $k = 2$}
				\addplot+ table[x=ndof,y=eta]{convhist_ex1_dG_adap_k_3_alpha_1000_tau_2.dat};
				\addlegendentry{\tiny $k = 3$}
				\addplot+ table[x=ndof,y=eta]{convhist_ex1_dG_unif_k_0_alpha_1_tau_2.dat};
				\addplot+ table[x=ndof,y=eta]{convhist_ex1_dG_unif_k_1_alpha_10_tau_2.dat};
				\addplot+ table[x=ndof,y=eta]{convhist_ex1_dG_unif_k_2_alpha_100_tau_2.dat};
				\addplot+ table[x=ndof,y=eta]{convhist_ex1_dG_unif_k_3_alpha_1000_tau_2.dat};
				\addplot+ [dotted,mark=none,black] coordinates{(1e1,6e-1) (1e10,6e-4)};
				\node(z) at  (axis cs:1e5,2e-2)
				[above] {\tiny $\mathcal O (\mathrm{ndof}^{-1/3})$};
				\addplot+ [dotted,mark=none,black] coordinates{(1e1,6e2) (1e6,6e-6)};
				\node(z) at  (axis cs:3e3,1e-2)
				[above] {\tiny $\mathcal O (\mathrm{ndof}^{-8/5})$};
			\end{loglogaxis}
		\end{tikzpicture}
	\end{subfigure}
	\begin{subfigure}[b]{0.48\textwidth}
		\begin{tikzpicture}[scale=0.78]
			\begin{loglogaxis}[legend pos=south west,legend cell align=left,
				legend style={fill=none},
				ymin=6e-5,ymax=1e0,
				xmin=20,xmax=1e6,
				ytick={1e-7,1e-6,1e-5,1e-4,1e-3,1e-2,1e-1,1e0,1e1,1e2},
				xtick={1e1,1e2,1e3,1e4,1e5,1e6,1e7}]
				\pgfplotsset{
					cycle list={%
						{Black, mark=diamond*, mark size=1.5pt},
						{Cyan, mark=*, mark size=1.5pt},
						{Orange, mark=square*, mark size=1.5pt},
						{Blue, mark=triangle*, mark size=1.5pt},
						{Black, dashed, mark=diamond*, mark size=1.5pt},
						{Cyan, dashed, mark=*, mark size=1.5pt},
						{Orange, dashed, mark=square*, mark size=1.5pt},
						{Blue, dashed, mark=triangle*, mark size=1.5pt},
						{Mulberry, dashed, mark=pentagon*, mark size=1.5pt}
					},
					legend style={
						at={(0.05,0.05)}, anchor=south west},
					font=\sffamily\scriptsize,
					xlabel=ndof,ylabel=
				}
				\addplot+ table[x=ndof,y=max_error]{convhist_ex1_dG_adap_k_0_alpha_1_tau_2.dat};
				\addlegendentry{\tiny $k = 0$}
				\addplot+ table[x=ndof,y=max_error]{convhist_ex1_dG_adap_k_1_alpha_10_tau_2.dat};
				\addlegendentry{\tiny $k = 1$}
				\addplot+ table[x=ndof,y=max_error]{convhist_ex1_dG_adap_k_2_alpha_100_tau_2.dat};
				\addlegendentry{\tiny $k = 2$}
				\addplot+ table[x=ndof,y=max_error]{convhist_ex1_dG_adap_k_3_alpha_1000_tau_2.dat};
				\addlegendentry{\tiny $k = 3$}
				\addplot+ table[x=ndof,y=max_error]{convhist_ex1_dG_unif_k_0_alpha_1_tau_2.dat};
				\addplot+ table[x=ndof,y=max_error]{convhist_ex1_dG_unif_k_1_alpha_10_tau_2.dat};
				\addplot+ table[x=ndof,y=max_error]{convhist_ex1_dG_unif_k_2_alpha_100_tau_2.dat};
				\addplot+ table[x=ndof,y=max_error]{convhist_ex1_dG_unif_k_3_alpha_1000_tau_2.dat};
				\addplot+ [dotted,mark=none,black] coordinates{(1e1,3e0) (1e7,3e-4)};
				\node(z) at  (axis cs:1e5,5e-3)
				[above] {\tiny $\mathcal O (\mathrm{ndof}^{-2/3})$};
				\addplot+ [dotted,mark=none,black] coordinates{(1e1,1e1) (1e6,1e-5)};
				\node(z) at  (axis cs:2e5,1e-4)
				[above] {\tiny $\mathcal O (\mathrm{ndof}^{-6/5})$};
				\addplot+ [dotted,mark=none,black] coordinates{(1e1,3e1) (1e6,3e-7)};
				\node(z) at  (axis cs:8e3,1e-4)
				[above] {\tiny $\mathcal O (\mathrm{ndof}^{-8/5})$};
			\end{loglogaxis}
		\end{tikzpicture}
	\end{subfigure}
	\caption{Convergence history of $\Phi_{\mathrm{nc}}(u_h)$ (left) and the $L^\infty$ error (right) for the NC-FEM in the second experiment with various $k$ and $\alpha = 10^{k}$ (solid lines = adaptive, dashed lines = uniform).
		\label{f:exp1-conv-dG}
	}
\end{figure}
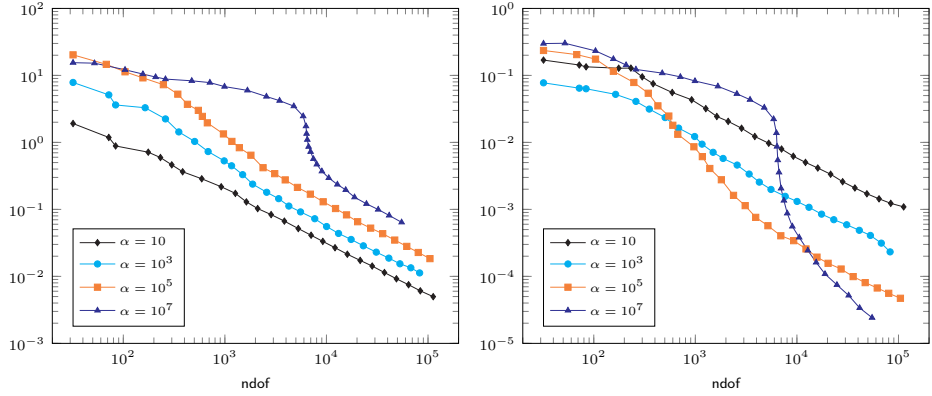
\begin{figure}
	\centering
	\begin{subfigure}[b]{0.48\textwidth}
		\begin{tikzpicture}[scale=0.78]
			\begin{loglogaxis}[legend pos=south west,legend cell align=left,
				legend style={fill=none},
				ymin=1e-3,ymax=1e2,
				xmin=20,xmax=2e5,
				ytick={1e-7,1e-6,1e-5,1e-4,1e-3,1e-2,1e-1,1e0,1e1,1e2,1e3},
				xtick={1e1,1e2,1e3,1e4,1e5,1e6,1e7}]
				\pgfplotsset{
					cycle list={%
						{Black, mark=diamond*, mark size=1.5pt},
						{Cyan, mark=*, mark size=1.5pt},
						{Orange, mark=square*, mark size=1.5pt},
						{Blue, mark=triangle*, mark size=1.5pt},
						{Black, dashed, mark=diamond*, mark size=1.5pt},
						{Cyan, dashed, mark=*, mark size=1.5pt},
						{Orange, dashed, mark=square*, mark size=1.5pt},
						{Blue, dashed, mark=triangle*, mark size=1.5pt},
						{Mulberry, dashed, mark=pentagon*, mark size=1.5pt}
					},
					legend style={
						at={(0.05,0.05)}, anchor=south west},
					font=\sffamily\scriptsize,
					xlabel=ndof,ylabel=
				}
				\addplot+ table[x=ndof,y=GUB]{convhist_ex1_bfs_adap_alpha_10_tau_2.dat};
				\addlegendentry{\tiny $\alpha = 10$}
				\addplot+ table[x=ndof,y=GUB]{convhist_ex1_bfs_adap_alpha_1000_tau_2.dat};
				\addlegendentry{\tiny $\alpha = 10^3$}
				\addplot+ table[x=ndof,y=GUB]{convhist_ex1_bfs_adap_alpha_100000_tau_2.dat};
				\addlegendentry{\tiny $\alpha = 10^5$}
				\addplot+ table[x=ndof,y=GUB]{convhist_ex1_bfs_adap_alpha_10000000_tau_2.dat};
				\addlegendentry{\tiny $\alpha = 10^7$}
			\end{loglogaxis}
		\end{tikzpicture}
	\end{subfigure}
	\begin{subfigure}[b]{0.48\textwidth}
		\begin{tikzpicture}[scale=0.78]
			\begin{loglogaxis}[legend pos=south west,legend cell align=left,
				legend style={fill=none},
				ymin=1e-5,ymax=1e0,
				xmin=20,xmax=2e5,
				ytick={1e-7,1e-6,1e-5,1e-4,1e-3,1e-2,1e-1,1e0,1e1,1e2},
				xtick={1e1,1e2,1e3,1e4,1e5,1e6,1e7}]
				\pgfplotsset{
					cycle list={%
						{Black, mark=diamond*, mark size=1.5pt},
						{Cyan, mark=*, mark size=1.5pt},
						{Orange, mark=square*, mark size=1.5pt},
						{Blue, mark=triangle*, mark size=1.5pt},
						{Black, dashed, mark=diamond*, mark size=1.5pt},
						{Cyan, dashed, mark=*, mark size=1.5pt},
						{Orange, dashed, mark=square*, mark size=1.5pt},
						{Blue, dashed, mark=triangle*, mark size=1.5pt},
						{Mulberry, dashed, mark=pentagon*, mark size=1.5pt}
					},
					legend style={
						at={(0.05,0.05)}, anchor=south west},
					font=\sffamily\scriptsize,
					xlabel=ndof,ylabel=
				}
				\addplot+ table[x=ndof,y=max_error]{convhist_ex1_bfs_adap_alpha_10_tau_2.dat};
				\addlegendentry{\tiny $\alpha = 10$}
				\addplot+ table[x=ndof,y=max_error]{convhist_ex1_bfs_adap_alpha_1000_tau_2.dat};
				\addlegendentry{\tiny $\alpha = 10^3$}
				\addplot+ table[x=ndof,y=max_error]{convhist_ex1_bfs_adap_alpha_100000_tau_2.dat};
				\addlegendentry{\tiny $\alpha = 10^5$}
				\addplot+ table[x=ndof,y=max_error]{convhist_ex1_bfs_adap_alpha_10000000_tau_2.dat};
				\addlegendentry{\tiny $\alpha = 10^7$}
			\end{loglogaxis}
		\end{tikzpicture}
	\end{subfigure}
	\caption{Convergence history of $\Phi(u_h)$ (left) and the $L^\infty$ error (right) for the BFS-FEM in the second experiment with various $\alpha$ on adaptive meshes.
		\label{f:exp1-conv-bfs-alpha}
	}
\end{figure}

\subsection{Second experiment}
In this benchmark,
we approximate the exact solution
\begin{align*}
	u(r,\varphi) \coloneqq r^{2/3} \sin(2\varphi/3)
\end{align*}
in polar coordinates to \eqref{def:PDE} in the L-shaped domain $\Omega = (-1,1)^2 \setminus ([0,1] \times [-1,0])$ with the coefficient
$$
A(r,\varphi) \coloneqq \begin{pmatrix}
	1 + 5r^{1/2} & r^{2}/2\\ r^{2}/2 & 1 + 5r^{1/2}
\end{pmatrix}
$$
in polar coordinates
and right-hand side $f(x) \coloneqq -A(x) : \D^2 u(x)$. The solution belongs to $H^{5/3-\delta}(\Omega)$ for any $\delta > 0$. This example appeared in \cite{QiuZhang2020}.
The parameter $\alpha$ is set to $10^k$, where $k$ denotes the order of the discretization ($k = 1$ for the BFS-FEM) for improved accuracy with higher polynomial degrees.
\Cref{f:exp1-conv-bfs}(a) displays the convergence history of the errors in standard norms and the GUB $\Phi(u_h)$ for the BFS-FEM.
Throughout this example, these errors and $\Phi(u_h)$ appear to converge with the same convergence rates.
Uniform mesh-refinements lead to the convergence rate 1/4 for these quantities.
The adaptive algorithm refines towards the reentrant corner as displayed in \Cref{f:exp1-conv-bfs}(b).
This leads to the improved convergence rate $4/5$ for the displayed quantities.
Thus, the (optimal) convergence rate $1$ in the smooth case was not recovered by adaptive computations. However, it is unclear whether this is possible.
This observation is consistent with the results from the NC-FEM displayed in \Cref{f:exp1-conv-dG}.
Higher polynomial degrees provide improved convergence rates.
We observed $1/3$ for $k = 0$, $2/3$ for $k = 1$, $6/5$ for $k = 2$, and $2$ for $k = 3$.
In \Cref{f:exp1-conv-bfs-alpha}, the influence of the parameter $\alpha$ on the convergence $\Phi(u_h)$ and the $L^\infty$ error for the BFS-FEM is studied.
While the convergence rates are comparable according to expectation, we observed more accurate discrete solutions (with respect to the $L^\infty$ error) for larger $\alpha$ after a preasymptotic regime. However, large $\alpha$ leads to larger GUB $\Phi(u_h)$.
In conclusion, it appears that an adaptive choice of $\alpha$ may improve the convergence rates of the $L^\infty$ error, but the coupling between $\alpha$ and the number of degrees of freedom remains unclear. Similar observations can be made in undisplayed results for the NC-FEM.

\subsection{Third experiment}
\begin{figure}
	\centering
	\begin{subfigure}[b]{0.48\textwidth}
		\begin{tikzpicture}[scale=0.78]
			\begin{loglogaxis}[legend pos=south west,legend cell align=left,
				legend style={fill=none},
				ymin=1e-4,ymax=1e-1,
				xmin=20,xmax=3e5,
				ytick={1e-7,1e-6,1e-5,1e-4,1e-3,1e-2,1e-1,1e0,1e1,1e2},
				xtick={1e1,1e2,1e3,1e4,1e5,1e6,1e7}]
				\pgfplotsset{
					cycle list={%
						{Black, mark=diamond*, mark size=1.5pt},
						{Black, dashed, mark=diamond*, mark size=1.5pt},
					},
					legend style={
						at={(0.05,0.05)}, anchor=south west},
					font=\sffamily\scriptsize,
					xlabel=ndof,ylabel=
				}
				\addplot+ table[x=ndof,y=GUB]{convhist_ex2_bfs_adap_alpha_10_tau_2.dat};
				\addlegendentry{\tiny $\Phi(u_h)$}
				\addplot+ table[x=ndof,y=GUB]{convhist_ex2_bfs_unif_alpha_10_tau_2.dat};
				\addplot+ [dotted,mark=none,black] coordinates{(1e1,5e-2) (1e11,5e-4)};
				\node(z) at  (axis cs:1e4,1e-2)
				[above] {\tiny $\mathcal O (\mathrm{ndof}^{-1/5})$};
				\addplot+ [dotted,mark=none,black] coordinates{(1e1,3e-1) (1e10,3e-7)};
				\node(z) at  (axis cs:6e4,1e-3)
				[above] {\tiny $\mathcal O (\mathrm{ndof}^{-2/3})$};
			\end{loglogaxis}
		\end{tikzpicture}
	\end{subfigure}
	\begin{subfigure}[b]{0.48\textwidth}
		\begin{tikzpicture}[scale=0.78]
			\begin{loglogaxis}[legend pos=south west,legend cell align=left,
				legend style={fill=none},
				ymin=1e-4,ymax=1e-1,
				xmin=20,xmax=1e6,
				ytick={1e-7,1e-6,1e-5,1e-4,1e-3,1e-2,1e-1,1e0,1e1,1e2,1e3},
				xtick={1e1,1e2,1e3,1e4,1e5,1e6,1e7}]
				\pgfplotsset{
					cycle list={%
						{Black, mark=diamond*, mark size=1.5pt},
						{Cyan, mark=*, mark size=1.5pt},
						{Orange, mark=square*, mark size=1.5pt},
						{Blue, mark=triangle*, mark size=1.5pt},
						{Black, dashed, mark=diamond*, mark size=1.5pt},
						{Cyan, dashed, mark=*, mark size=1.5pt},
						{Orange, dashed, mark=square*, mark size=1.5pt},
						{Blue, dashed, mark=triangle*, mark size=1.5pt},
						{Mulberry, dashed, mark=pentagon*, mark size=1.5pt}
					},
					legend style={
						at={(0.05,0.05)}, anchor=south west},
					font=\sffamily\scriptsize,
					xlabel=ndof,ylabel=
				}
				\addplot+ table[x=ndof,y=eta]{convhist_ex2_dG_adap_k_0_alpha_10_tau_2.dat};
				\addlegendentry{\tiny $k = 0$}
				\addplot+ table[x=ndof,y=eta]{convhist_ex2_dG_adap_k_1_alpha_10_tau_2.dat};
				\addlegendentry{\tiny $k = 1$}
				\addplot+ table[x=ndof,y=eta]{convhist_ex2_dG_adap_k_2_alpha_10_tau_2.dat};
				\addlegendentry{\tiny $k = 2$}
				\addplot+ table[x=ndof,y=eta]{convhist_ex2_dG_adap_k_3_alpha_10_tau_2.dat};
				\addlegendentry{\tiny $k = 3$}
				\addplot+ table[x=ndof,y=eta]{convhist_ex2_dG_unif_k_0_alpha_10_tau_2.dat};
				\addplot+ table[x=ndof,y=eta]{convhist_ex2_dG_unif_k_1_alpha_10_tau_2.dat};
				\addplot+ table[x=ndof,y=eta]{convhist_ex2_dG_unif_k_2_alpha_10_tau_2.dat};
				\addplot+ table[x=ndof,y=eta]{convhist_ex2_dG_unif_k_3_alpha_10_tau_2.dat};
				\addplot+ [dotted,mark=none,black] coordinates{(1e1,2e-1) (1e10,2e-4)};
				\node(z) at  (axis cs:1e5,7e-3)
				[above] {\tiny $\mathcal O (\mathrm{ndof}^{-1/3})$};
				\addplot+ [dotted,mark=none,black] coordinates{(1e1,6e-1) (1e10,6e-7)};
				\node(z) at  (axis cs:1.5e5,1e-3)
				[above] {\tiny $\mathcal O (\mathrm{ndof}^{-2/3})$};
				\addplot+ [dotted,mark=none,black] coordinates{(1e1,1e1) (1e6,1e-7)};
				\node(z) at  (axis cs:3e3,7e-4)
				[above] {\tiny $\mathcal O (\mathrm{ndof}^{-8/5})$};
			\end{loglogaxis}
		\end{tikzpicture}
	\end{subfigure}
	\caption{Convergence history of $\Phi(u_h)$ for the BFS-FEM (left) and $\Phi_{\mathrm{nc}}(u_h)$ for the NC-FEM (right) in the third experiment with $\alpha \coloneqq 10$ (solid lines = adaptive, dashed lines = uniform).}
	\label{f:exp2-conv}
\end{figure}
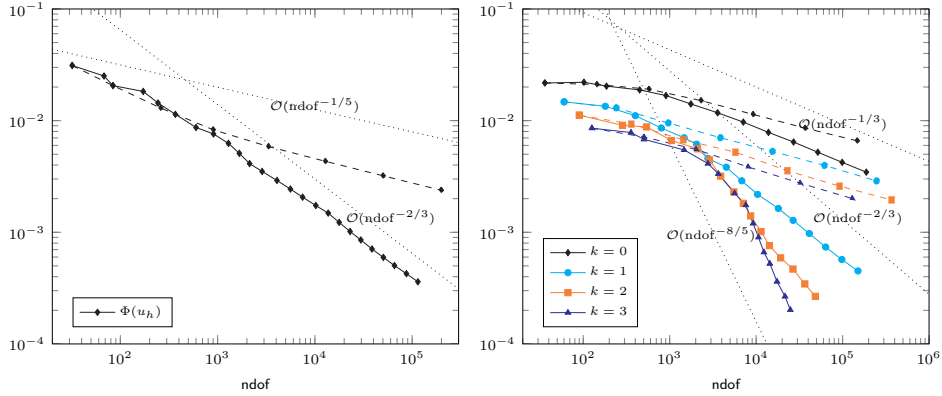

In this benchmark, we approximate the unknown solution to \eqref{def:PDE} in the L-shaped domain $\Omega = (-1,1)^2 \setminus ([0,1] \times [-1,0])$ with the coefficient
\begin{align*}
	A(r,\varphi) \coloneqq \begin{pmatrix}
		15 - r^{1/2} & 1\\ 1 & 3 + r.^{1/2}
	\end{pmatrix}
\end{align*}
in polar coordinates,
right-hand side $f \equiv 1$, and homogenous Dirichlet data $g \equiv 0$.
Conforming methods can provide unconditional information on the $L^\infty$ error by the evaluation of the GUB $\Phi(u_h)$.
Due to the lack of an exact solution, the quantities of interest are $\Phi(u_h)$ for the BFS-FEM and $\Phi_{\mathrm{nc}}(u_h)$ for the NC-FEM. In this example, we set $\alpha \coloneqq 10$.
\Cref{f:exp2-conv} displays the convergence history of the aforementioned quantities with similar results as in the previous example.
Uniform mesh refinements lead to the convergence rate $1/5$, while adaptive computation refines towards the reentrant corner and provide improved convergence rates. For the NC-FEM, we observe $1/3$ for $k = 0$, $2/3$ for $k = 1$ (as well as for the BFS-FEM), around $1$ for $k = 2$, and $8/5$ for $k = 3$.

\begin{figure}
	\centering
	\begin{subfigure}[b]{0.48\textwidth}
		\begin{tikzpicture}[scale=0.78]
			\begin{loglogaxis}[legend pos=south west,legend cell align=left,
				legend style={fill=none},
				ymin=6e-3,ymax=1e1,
				xmin=20,xmax=1e5,
				ytick={1e-7,1e-6,1e-5,1e-4,1e-3,1e-2,1e-1,1e0,1e1,1e2},
				xtick={1e1,1e2,1e3,1e4,1e5,1e6,1e7}]
				\pgfplotsset{
					cycle list={%
						{Black, mark=diamond*, mark size=1.5pt},
						{Cyan, mark=*, mark size=1.5pt},
						{Orange, mark=square*, mark size=1.5pt},
						{Blue, mark=triangle*, mark size=1.5pt},
						{Black, dashed, mark=diamond*, mark size=1.5pt},
						{Cyan, dashed, mark=*, mark size=1.5pt},
						{Orange, dashed, mark=square*, mark size=1.5pt},
						{Blue, dashed, mark=triangle*, mark size=1.5pt},
						{Mulberry, dashed, mark=pentagon*, mark size=1.5pt}
					},
					legend style={
						at={(0.05,0.05)}, anchor=south west},
					font=\sffamily\scriptsize,
					xlabel=ndof,ylabel=
				}
				\addplot+ table[x=ndof,y=GUB]{convhist_ex3_bfs_adap_alpha_10_tau_2.dat};
				\addlegendentry{\tiny $\Phi(u_h)$}
				\addplot+ 	table[x=ndof,y=max_error]{convhist_ex3_bfs_adap_alpha_10_tau_2.dat};
				\addlegendentry{\tiny $L^\infty$ error}
				\addplot+ 	table[x=ndof,y=L2_error]{convhist_ex3_bfs_adap_alpha_10_tau_2.dat};
				\addlegendentry{\tiny $L^2$ error}
				\addplot+ 	table[x=ndof,y=H1_error]{convhist_ex3_bfs_adap_alpha_10_tau_2.dat};
				\addlegendentry{\tiny $H^1$ error}
				\addplot+ table[x=ndof,y=GUB]{convhist_ex3_bfs_unif_alpha_10_tau_2.dat};
				\addplot+ 	table[x=ndof,y=max_error]{convhist_ex3_bfs_unif_alpha_10_tau_2.dat};
				\addplot+ 	table[x=ndof,y=L2_error]{convhist_ex3_bfs_unif_alpha_10_tau_2.dat};
				\addplot+ 	table[x=ndof,y=H1_error]{convhist_ex3_bfs_unif_alpha_10_tau_2.dat};
				\addplot+ [dotted,mark=none,black] coordinates{(1e1,2e0) (1e13,2e-2)};
				\node(z) at  (axis cs:1e4,4.5e-1)
				[above] {\tiny $\mathcal O (\mathrm{ndof}^{-1/6})$};
				\addplot+ [dotted,mark=none,black] coordinates{(1e1,2e0) (1e10,2e-6)};
				\node(z) at  (axis cs:1e4,1e-2)
				[above] {\tiny $\mathcal O (\mathrm{ndof}^{-2/3})$};
			\end{loglogaxis}
		\end{tikzpicture}
	\end{subfigure}
	\begin{subfigure}[b]{0.48\textwidth}
		\centering
		\includegraphics[scale=0.57]{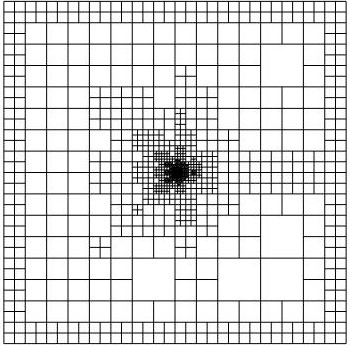}
	\end{subfigure}
	\caption{Convergence history (left) and adaptive mesh with 2122 elements (right) of the BFS-FEM for the fourth experiment with $\alpha = 10$ (solid lines = adaptive, dashed lines = uniform).}
	\label{f:exp3-conv-bfs}
\end{figure}
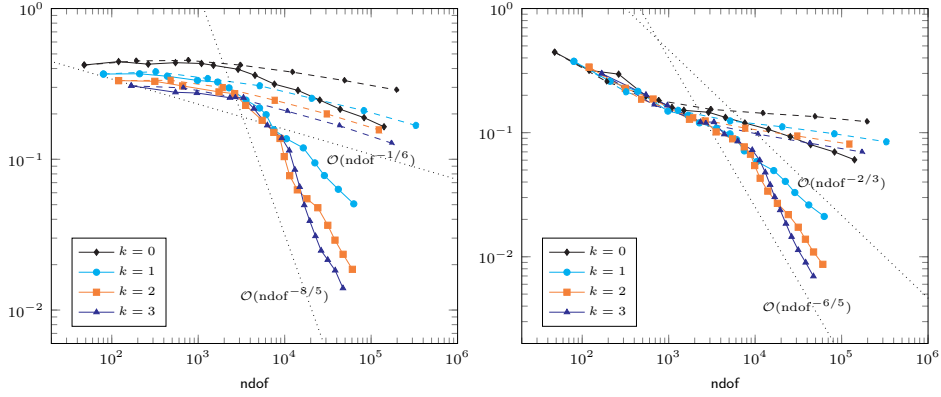
\begin{figure}
	\centering
	\begin{subfigure}[b]{0.48\textwidth}
		\begin{tikzpicture}[scale=0.78]
			\begin{loglogaxis}[legend pos=south west,legend cell align=left,
				legend style={fill=none},
				ymin=6e-3,ymax=1e0,
				xmin=20,xmax=1e6,
				ytick={1e-7,1e-6,1e-5,1e-4,1e-3,1e-2,1e-1,1e0,1e1,1e2,1e3},
				xtick={1e1,1e2,1e3,1e4,1e5,1e6,1e7}]
				\pgfplotsset{
					cycle list={%
						{Black, mark=diamond*, mark size=1.5pt},
						{Cyan, mark=*, mark size=1.5pt},
						{Orange, mark=square*, mark size=1.5pt},
						{Blue, mark=triangle*, mark size=1.5pt},
						{Black, dashed, mark=diamond*, mark size=1.5pt},
						{Cyan, dashed, mark=*, mark size=1.5pt},
						{Orange, dashed, mark=square*, mark size=1.5pt},
						{Blue, dashed, mark=triangle*, mark size=1.5pt},
						{Mulberry, dashed, mark=pentagon*, mark size=1.5pt}
					},
					legend style={
						at={(0.05,0.05)}, anchor=south west},
					font=\sffamily\scriptsize,
					xlabel=ndof,ylabel=
				}
				\addplot+ table[x=ndof,y=eta]{convhist_ex3_dG_adap_k_0_alpha_10_tau_2.dat};
				\addlegendentry{\tiny $k = 0$}
				\addplot+ table[x=ndof,y=eta]{convhist_ex3_dG_adap_k_1_alpha_10_tau_2.dat};
				\addlegendentry{\tiny $k = 1$}
				\addplot+ table[x=ndof,y=eta]{convhist_ex3_dG_adap_k_2_alpha_10_tau_2.dat};
				\addlegendentry{\tiny $k = 2$}
				\addplot+ table[x=ndof,y=eta]{convhist_ex3_dG_adap_k_3_alpha_10_tau_2.dat};
				\addlegendentry{\tiny $k = 3$}
				\addplot+ table[x=ndof,y=eta]{convhist_ex3_dG_unif_k_0_alpha_10_tau_2.dat};
				\addplot+ table[x=ndof,y=eta]{convhist_ex3_dG_unif_k_1_alpha_10_tau_2.dat};
				\addplot+ table[x=ndof,y=eta]{convhist_ex3_dG_unif_k_2_alpha_10_tau_2.dat};
				\addplot+ table[x=ndof,y=eta]{convhist_ex3_dG_unif_k_3_alpha_10_tau_2.dat};
				\addplot+ [dotted,mark=none,black] coordinates{(1e1,5e-1) (1e13,5e-3)};
				\node(z) at  (axis cs:1e5,8e-2)
				[above] {\tiny $\mathcal O (\mathrm{ndof}^{-1/6})$};
				\addplot+ [dotted,mark=none,black] coordinates{(1e1,2e3) (1e6,2e-5)};
				\node(z) at  (axis cs:1e4,1e-2)
				[above] {\tiny $\mathcal O (\mathrm{ndof}^{-8/5})$};
			\end{loglogaxis}
		\end{tikzpicture}
	\end{subfigure}
	\begin{subfigure}[b]{0.48\textwidth}
		\begin{tikzpicture}[scale=0.78]
			\begin{loglogaxis}[legend pos=south west,legend cell align=left,
				legend style={fill=none},
				ymin=2e-3,ymax=1e0,
				xmin=20,xmax=1e6,
				ytick={1e-7,1e-6,1e-5,1e-4,1e-3,1e-2,1e-1,1e0,1e1,1e2},
				xtick={1e1,1e2,1e3,1e4,1e5,1e6,1e7}]
				\pgfplotsset{
					cycle list={%
						{Black, mark=diamond*, mark size=1.5pt},
						{Cyan, mark=*, mark size=1.5pt},
						{Orange, mark=square*, mark size=1.5pt},
						{Blue, mark=triangle*, mark size=1.5pt},
						{Black, dashed, mark=diamond*, mark size=1.5pt},
						{Cyan, dashed, mark=*, mark size=1.5pt},
						{Orange, dashed, mark=square*, mark size=1.5pt},
						{Blue, dashed, mark=triangle*, mark size=1.5pt},
						{Mulberry, dashed, mark=pentagon*, mark size=1.5pt}
					},
					legend style={
						at={(0.05,0.05)}, anchor=south west},
					font=\sffamily\scriptsize,
					xlabel=ndof,ylabel=
				}
				\addplot+ table[x=ndof,y=max_error]{convhist_ex3_dG_adap_k_0_alpha_10_tau_2.dat};
				\addlegendentry{\tiny $k = 0$}
				\addplot+ table[x=ndof,y=max_error]{convhist_ex3_dG_adap_k_1_alpha_10_tau_2.dat};
				\addlegendentry{\tiny $k = 1$}
				\addplot+ table[x=ndof,y=max_error]{convhist_ex3_dG_adap_k_2_alpha_10_tau_2.dat};
				\addlegendentry{\tiny $k = 2$}
				\addplot+ table[x=ndof,y=max_error]{convhist_ex3_dG_adap_k_3_alpha_10_tau_2.dat};
				\addlegendentry{\tiny $k = 3$}
				\addplot+ table[x=ndof,y=max_error]{convhist_ex3_dG_unif_k_0_alpha_10_tau_2.dat};
				\addplot+ table[x=ndof,y=max_error]{convhist_ex3_dG_unif_k_1_alpha_10_tau_2.dat};
				\addplot+ table[x=ndof,y=max_error]{convhist_ex3_dG_unif_k_2_alpha_10_tau_2.dat};
				\addplot+ table[x=ndof,y=max_error]{convhist_ex3_dG_unif_k_3_alpha_10_tau_2.dat};
				\addplot+ [dotted,mark=none,black] coordinates{(1e1,1e1) (1e7,1e-3)};
				\node(z) at  (axis cs:1e5,3e-2)
				[above] {\tiny $\mathcal O (\mathrm{ndof}^{-2/3})$};
				\addplot+ [dotted,mark=none,black] coordinates{(1e1,1e2) (1e6,1e-4)};
				\node(z) at  (axis cs:4e4,3e-3)
				[above] {\tiny $\mathcal O (\mathrm{ndof}^{-6/5})$};
			\end{loglogaxis}
		\end{tikzpicture}
	\end{subfigure}
	\caption{Convergence history of $\Phi_{\mathrm{nc}}(u_h)$ (left) and the $L^\infty$ error (right) for the NC-FEM in the fourth experiment with various $k$ and $\alpha = 10$ (solid lines = adaptive, dashed lines = uniform).
		\label{f:exp3-conv-dG}
	}
\end{figure}

\subsection{Fourth experiment}
In this benchmark, we approximate the exact solution
\begin{align*}
	u(r,\varphi) \coloneqq r^{1/2}\sin(\varphi/2) - r^2\sin^2(\varphi)
\end{align*}
to \eqref{def:PDE} in the slit domain $\Omega \coloneqq (-1,1)^2 \setminus ([0,1] \times \{0\})$ with the discontinuous coefficient, for $x = (x_1,x_2)$,
\begin{align*}
	A(x) \coloneqq \begin{cases}
		\begin{pmatrix}
			1 + 5|x|^{1/2} & |x|^2/2\\ |x|^2/2 & 1 + 5|x|^{1/2}
		\end{pmatrix} 
		&\mbox{if } x_1 \leq x_2,\\
		(1 + |x - (0,-1)|^{1/3})\mathrm{I}_2 &\mbox{otherwise}
	\end{cases}
\end{align*}
and right-hand side $f(x) = 1$ if $x_1 \geq x_2$ and $f(x) = (1 + |x-(-1,1)|^{1/3})$ otherwise.
The function $u$ belongs to $H^{3/2-\delta}(\Omega)$ for any $\delta > 0$.
The convergence analysis of this paper does not apply to this example because $A$ is discontinuous and $\Omega$ is not a Lipschitz domain.
Nevertheless, the ABP maximum principle from \Cref{thm:ABP-maximum-principle} applies to this example as well, so $\Phi(u_h)$ is a guaranteed upper bound for $\|u - u_h\|_{L^\infty(\Omega)}$ provided a discrete function $u_h \in W^{2,n}(\Omega)$ is given.
Notice that this requires the information $u \in C(\overline{\Omega}) \cap H^2_\mathrm{loc}(\Omega)$ on the exact solution $u$.
In this example, we set $\alpha = 10$.
The results for the BFS-FEM and NC-FEM displayed in \Cref{f:exp3-conv-bfs}--\ref{f:exp3-conv-dG} match the observations of the previous experiments, although this example is not covered by the theory.
\Cref{f:exp3-conv-bfs}(b) shows that the adaptive algorithm refines towards the reentrant corner, but not along the set of discontinuity of the coefficient $A$, which is the diagonal from the bottom left to the upper right corner.

\subsection{Conclusion}
In all computer experiments, we observed that $\Phi(u_h)$ is indeed a guaranteed bound for the $L^\infty$ error if $u_h$ is a conforming approximation. For nonsmooth exact solutions $u$, $\Phi(u_h)$ appears to be efficient as well. In this case, the efficiency index $\Phi(u_h)/\|u - u_h\|_{L^\infty(\Omega)}$ depends on the parameter $\alpha$, where larger $\alpha$ leads to a larger index. A similar behavior is observed for the a~posteriori error estimator $\Phi_{\mathrm{nc}}(u_h)$ for the NC-FEM.
Adaptive computations lead to improved convergence rates of the minimizing functional and provides significant improvements to the convergence of $u_h$ towards $u$ in comparison to uniform mesh-refinements.
However, the convergence rates in the smooth case could not be recovered. It remains unclear whether this is possible or how to achieve it.
A straightforward adaptive strategy is unavailable because a part of the a posteriori error control cannot be localize.
We note that, for discretizations of higher polynomial order, the $L^n$ contribution in the objective functionals will dominate the boundary residual. Thus, it is expected that a localization of the $L^n$ contribution is sufficient to drive the adaptive mesh-refining algorithm for $k >> 1$.
\printbibliography

@article {FengJensen2017,
    AUTHOR = {Feng, Xiaobing and Jensen, Max},
     TITLE = {Convergent semi-{L}agrangian methods for the {M}onge-{A}mp\`ere
              equation on unstructured grids},
   JOURNAL = {SIAM J. Numer. Anal.},
  FJOURNAL = {SIAM Journal on Numerical Analysis},
    VOLUME = {55},
      YEAR = {2017},
    NUMBER = {2},
     PAGES = {691--712},
      ISSN = {0036-1429},
   MRCLASS = {65N06 (35D40 35J25 35J60 35J96 65N12)},
  MRNUMBER = {3623696},
MRREVIEWER = {Thomas Lee Lewis},
       DOI = {10.1137/16M1061709},
       URL = {https://doi.org/10.1137/16M1061709},
}

@article {BarlesSouganidis1991,
	AUTHOR = {Barles, G. and Souganidis, P. E.},
	TITLE = {Convergence of approximation schemes for fully nonlinear
	second order equations},
	JOURNAL = {Asymptotic Anal.},
	FJOURNAL = {Asymptotic Analysis},
	VOLUME = {4},
	YEAR = {1991},
	NUMBER = {3},
	PAGES = {271--283},
	ISSN = {0921-7134},
	MRCLASS = {35K55 (35A40 35J60 65M12)},
	MRNUMBER = {1115933},
}

@article {MotzkinWasow1953,
	AUTHOR = {Motzkin, T. S. and Wasow, W.},
	TITLE = {On the approximation of linear elliptic differential equations
	by difference equations with positive coefficients},
	JOURNAL = {J. Math. Physics},
	VOLUME = {31},
	YEAR = {1953},
	PAGES = {253--259},
	MRCLASS = {65.0X},
	MRNUMBER = {0052895},
	MRREVIEWER = {L. Bers},
}

@article {SmearsSueli2013,
	AUTHOR = {Smears, Iain and S\"{u}li, Endre},
	TITLE = {Discontinuous {G}alerkin finite element approximation of
	nondivergence form elliptic equations with {C}ord\`es
	coefficients},
	JOURNAL = {SIAM J. Numer. Anal.},
	FJOURNAL = {SIAM Journal on Numerical Analysis},
	VOLUME = {51},
	YEAR = {2013},
	NUMBER = {4},
	PAGES = {2088--2106},
	ISSN = {0036-1429},
	MRCLASS = {65N30 (35J25 65N12)},
	MRNUMBER = {3077903},
	MRREVIEWER = {Andreas Schr\"{o}der},
	DOI = {10.1137/120899613},
	URL = {https://doi.org/10.1137/120899613},
}

@article {SmearsSueli2014,
	AUTHOR = {Smears, Iain and S\"{u}li, Endre},
	TITLE = {Discontinuous {G}alerkin finite element approximation of
	{H}amilton-{J}acobi-{B}ellman equations with {C}ordes
	coefficients},
	JOURNAL = {SIAM J. Numer. Anal.},
	FJOURNAL = {SIAM Journal on Numerical Analysis},
	VOLUME = {52},
	YEAR = {2014},
	NUMBER = {2},
	PAGES = {993--1016},
	ISSN = {0036-1429},
	MRCLASS = {65N30 (35D35 35J60 47J25 49M25)},
	MRNUMBER = {3196952},
	MRREVIEWER = {Marius Ghergu},
	DOI = {10.1137/130909536},
	URL = {https://doi.org/10.1137/130909536},
}

@article {GallistlSueli2019,
	AUTHOR = {Gallistl, Dietmar and S\"{u}li, Endre},
	TITLE = {Mixed finite element approximation of the
	{H}amilton-{J}acobi-{B}ellman equation with {C}ordes
	coefficients},
	JOURNAL = {SIAM J. Numer. Anal.},
	FJOURNAL = {SIAM Journal on Numerical Analysis},
	VOLUME = {57},
	YEAR = {2019},
	NUMBER = {2},
	PAGES = {592--614},
	ISSN = {0036-1429},
	MRCLASS = {65N30 (65N12 65N15 65N50)},
	MRNUMBER = {3924618},
	MRREVIEWER = {Qingfang Liu},
	DOI = {10.1137/18M1192299},
	URL = {https://doi.org/10.1137/18M1192299},
}

@article {CaffarelliCrandallKocanSwiech1996,
    AUTHOR = {Caffarelli, L. and Crandall, M. G. and Kocan, M. and \'Swi\k{e}ch,
              A.},
     TITLE = {On viscosity solutions of fully nonlinear equations with
              measurable ingredients},
   JOURNAL = {Comm. Pure Appl. Math.},
  FJOURNAL = {Communications on Pure and Applied Mathematics},
    VOLUME = {49},
      YEAR = {1996},
    NUMBER = {4},
     PAGES = {365--397},
      ISSN = {0010-3640},
   MRCLASS = {35J60 (35D05 35D10)},
  MRNUMBER = {1376656},
MRREVIEWER = {Katsuyuki Ishii},
       DOI = {10.1002/(SICI)1097-0312(199604)49:4<365::AID-CPA3>3.3.CO;2-V},
       URL =
              {https://doi.org/10.1002/(SICI)1097-0312(199604)49:4<365::AID-CPA3>3.3.CO;2-V},
}

@book {GilbargTrudinger2001,
    AUTHOR = {Gilbarg, David and Trudinger, Neil S.},
     TITLE = {Elliptic partial differential equations of second order},
    SERIES = {Classics in Mathematics},
      NOTE = {Reprint of the 1998 edition},
 PUBLISHER = {Springer-Verlag, Berlin},
      YEAR = {2001},
     PAGES = {xiv+517},
      ISBN = {3-540-41160-7},
   MRCLASS = {35-02 (35Jxx)},
  MRNUMBER = {1814364},
}

@article {KoikeSwiech2009,
    AUTHOR = {Koike, Shigeaki and \'{S}wi\c{e}ch, Andrzej},
     TITLE = {Weak {H}arnack inequality for fully nonlinear uniformly
              elliptic {PDE} with unbounded ingredients},
   JOURNAL = {J. Math. Soc. Japan},
  FJOURNAL = {Journal of the Mathematical Society of Japan},
    VOLUME = {61},
      YEAR = {2009},
    NUMBER = {3},
     PAGES = {723--755},
      ISSN = {0025-5645},
   MRCLASS = {35J60 (35B45 49L25)},
  MRNUMBER = {2552914},
MRREVIEWER = {Niko M. Marola},
       URL = {http://projecteuclid.org/euclid.jmsj/1248961477},
}

@book {AdamsFournier2003,
	AUTHOR = {Adams, R. A. and Fournier, J. J. F.},
	TITLE = {Sobolev spaces},
	SERIES = {Pure and Applied Mathematics (Amsterdam)},
	VOLUME = {140},
	EDITION = {Second},
	PUBLISHER = {Elsevier/Academic Press, Amsterdam},
	YEAR = {2003},
	PAGES = {xiv+305},
	ISBN = {0-12-044143-8},
	MRCLASS = {46E35 (46-01 46-02 46B70 46Exx)},
	MRNUMBER = {2424078},
}

@article {KaweckiSmears2022,
    AUTHOR = {Kawecki, Ellya L. and Smears, Iain},
     TITLE = {Convergence of adaptive discontinuous {G}alerkin and
              {$C^0$}-interior penalty finite element methods for
              {H}amilton-{J}acobi-{B}ellman and {I}saacs equations},
   JOURNAL = {Found. Comput. Math.},
  FJOURNAL = {Foundations of Computational Mathematics. The Journal of the
              Society for the Foundations of Computational Mathematics},
    VOLUME = {22},
      YEAR = {2022},
    NUMBER = {2},
     PAGES = {315--364},
      ISSN = {1615-3375},
   MRCLASS = {65N15 (65N30 65N50)},
  MRNUMBER = {4407745},
MRREVIEWER = {Mahboub Baccouch},
       DOI = {10.1007/s10208-021-09493-0},
       URL = {https://doi.org/10.1007/s10208-021-09493-0},
}

@article {Alexandrov1966,
    AUTHOR = {Alexandrov, A. D.},
     TITLE = {The impossibility of general estimates for solutions and
			  of uniqueness conditions for linear
			  equations with norms weaker than in $L_n$},
   JOURNAL = {Vestnik Leningrad Univ.},
    VOLUME = {21},
      YEAR = {1966},
    NUMBER = {12},
     PAGES = {5--10}
}

@book {Grisvard2011,
    AUTHOR = {Grisvard, Pierre},
     TITLE = {Elliptic problems in nonsmooth domains},
    SERIES = {Classics in Applied Mathematics},
    VOLUME = {69},
 PUBLISHER = {SIAM, Philadelphia, PA},
      YEAR = {2011},
     PAGES = {xx+410},
      ISBN = {978-1-611972-02-3},
   MRCLASS = {35J25 (01A75 35-02)},
  MRNUMBER = {3396210},
       DOI = {10.1137/1.9781611972030.ch1},
       URL = {https://doi.org/10.1137/1.9781611972030.ch1},
}

@article {FengHenningsNeilan2017,
    AUTHOR = {Feng, Xiaobing and Hennings, Lauren and Neilan, Michael},
     TITLE = {Finite element methods for second order linear elliptic
              partial differential equations in non-divergence form},
   JOURNAL = {Math. Comp.},
  FJOURNAL = {Mathematics of Computation},
    VOLUME = {86},
      YEAR = {2017},
    NUMBER = {307},
     PAGES = {2025--2051},
      ISSN = {0025-5718},
   MRCLASS = {65N30 (35J25 65N12)},
  MRNUMBER = {3647950},
MRREVIEWER = {Alexandre Ern},
       DOI = {10.1090/mcom/3168},
       URL = {https://doi.org/10.1090/mcom/3168},
}

@article {BrennerGudiSung2010,
    AUTHOR = {Brenner, Susanne C. and Gudi, Thirupathi and Sung, Li-yeng},
     TITLE = {An a posteriori error estimator for a quadratic
              {$C^0$}-interior penalty method for the biharmonic problem},
   JOURNAL = {IMA J. Numer. Anal.},
  FJOURNAL = {IMA Journal of Numerical Analysis},
    VOLUME = {30},
      YEAR = {2010},
    NUMBER = {3},
     PAGES = {777--798},
      ISSN = {0272-4979},
   MRCLASS = {65N30 (35J40 65N15)},
  MRNUMBER = {2670114},
MRREVIEWER = {Jeffrey S. Ovall},
       DOI = {10.1093/imanum/drn057},
       URL = {https://doi.org/10.1093/imanum/drn057},
}

@article {Talenti1965,
    AUTHOR = {Talenti, Giorgio},
     TITLE = {Sopra una classe di equazioni ellittiche a coefficienti
              misurabili},
   JOURNAL = {Ann. Mat. Pura Appl. (4)},
  FJOURNAL = {Annali di Matematica Pura ed Applicata. Serie Quarta},
    VOLUME = {69},
      YEAR = {1965},
     PAGES = {285--304},
      ISSN = {0003-4622},
   MRCLASS = {35.42},
  MRNUMBER = {201816},
MRREVIEWER = {G. Geymonat},
       DOI = {10.1007/BF02414375},
       URL = {https://doi.org/10.1007/BF02414375},
}

@article {Nadirashvili1997,
    AUTHOR = {Nadirashvili, Nikolai},
     TITLE = {Nonuniqueness in the martingale problem and the {D}irichlet
              problem for uniformly elliptic operators},
   JOURNAL = {Ann. Scuola Norm. Sup. Pisa Cl. Sci. (4)},
  FJOURNAL = {Annali della Scuola Normale Superiore di Pisa. Classe di
              Scienze. Serie IV},
    VOLUME = {24},
      YEAR = {1997},
    NUMBER = {3},
     PAGES = {537--549},
      ISSN = {0391-173X},
   MRCLASS = {35J25 (60J60)},
  MRNUMBER = {1612401},
MRREVIEWER = {Joerg-Uwe Loebus},
       URL = {http://www.numdam.org/item?id=ASNSP_1997_4_24_3_537_0},
}

@article {Safonov1999,
    AUTHOR = {Safonov, Mikhail V.},
     TITLE = {Nonuniqueness for second-order elliptic equations with
              measurable coefficients},
   JOURNAL = {SIAM J. Math. Anal.},
  FJOURNAL = {SIAM Journal on Mathematical Analysis},
    VOLUME = {30},
      YEAR = {1999},
    NUMBER = {4},
     PAGES = {879--895},
      ISSN = {0036-1410},
   MRCLASS = {35J15 (35B20 35R05 60G44 60J60)},
  MRNUMBER = {1684729},
       DOI = {10.1137/S0036141096309046},
       URL = {https://doi.org/10.1137/S0036141096309046},
}

@article {Doktor1976,
    AUTHOR = {Doktor, Pavel},
     TITLE = {Approximation of domains with {L}ipschitzian boundary},
   JOURNAL = {\v{C}asopis P\v{e}st. Mat.},
  FJOURNAL = {\v{C}eskoslovensk\'{a} Akademie V\v{e}d. \v{C}asopis Pro P\v{e}stov\'{a}n\'{\i} Matematiky},
    VOLUME = {101},
      YEAR = {1976},
    NUMBER = {3},
     PAGES = {237--255},
      ISSN = {0528-2195},
   MRCLASS = {46E35},
  MRNUMBER = {0461122},
MRREVIEWER = {R. A. Adams},
}

@book {Ciarlet1978,
    AUTHOR = {Ciarlet, Philippe G.},
     TITLE = {The Finite Element Method for Elliptic Problems},
    SERIES = {Studies in Mathematics and its Applications},
    VOLUME = {4},
 PUBLISHER = {North-Holland},
   ADDRESS = {Amsterdam},
      YEAR = {1978},
}

@article {WorseyFarin1987,
    AUTHOR = {Worsey, A. J. and Farin, G.},
     TITLE = {An {$n$}-dimensional {C}lough-{T}ocher interpolant},
   JOURNAL = {Constr. Approx.},
  FJOURNAL = {Constructive Approximation. An International Journal for
              Approximations and Expansions},
    VOLUME = {3},
      YEAR = {1987},
    NUMBER = {2},
     PAGES = {99--110},
      ISSN = {0176-4276},
   MRCLASS = {41A05 (41A15 41A63)},
  MRNUMBER = {889547},
MRREVIEWER = {G. Baszenski},
       DOI = {10.1007/BF01890556},
       URL = {https://doi.org/10.1007/BF01890556},
}

@article {NochettoZhang2018,
    AUTHOR = {Nochetto, Ricardo H. and Zhang, Wujun},
     TITLE = {Discrete {ABP} estimate and convergence rates for linear
              elliptic equations in non-divergence form},
   JOURNAL = {Found. Comput. Math.},
  FJOURNAL = {Foundations of Computational Mathematics. The Journal of the
              Society for the Foundations of Computational Mathematics},
    VOLUME = {18},
      YEAR = {2018},
    NUMBER = {3},
     PAGES = {537--593},
      ISSN = {1615-3375},
   MRCLASS = {65N30 (35B50 35D35 35J25 65N15)},
  MRNUMBER = {3807356},
MRREVIEWER = {Mir Sajjad Hashemi},
       DOI = {10.1007/s10208-017-9347-y},
       URL = {https://doi.org/10.1007/s10208-017-9347-y},
}

@article {CloughTocher1965,
    AUTHOR = {Clough, R. W. and Tocher, J. L.},
     TITLE = {Finite Element Stiffness Matrices for Analysis of Plates in Bending},
   JOURNAL = {Proceedings of the Conference on Matrix Methods in Structural Mechanics},
      YEAR = {1965},
     PAGES = {515--545},
}

@article {DouglasDupontPercellScott1979,
    AUTHOR = {Douglas, Jr., Jim and Dupont, Todd and Percell, Peter and
              Scott, Ridgway},
     TITLE = {A family of {$C^{1}$} finite elements with optimal
              approximation properties for various {G}alerkin methods for
              2nd and 4th order problems},
   JOURNAL = {RAIRO Anal. Num\'{e}r.},
  FJOURNAL = {RAIRO Analyse Num\'{e}rique},
    VOLUME = {13},
      YEAR = {1979},
    NUMBER = {3},
     PAGES = {227--255},
      ISSN = {0399-0516},
   MRCLASS = {65N30},
  MRNUMBER = {543934},
MRREVIEWER = {Reinhard Scholz},
       DOI = {10.1051/m2an/1979130302271},
       URL = {https://doi.org/10.1051/m2an/1979130302271},
}

@article {GeorgoulisHoustonVirtanen2011,
    AUTHOR = {Georgoulis, Emmanuil H. and Houston, Paul and Virtanen, Juha},
     TITLE = {An {\it a posteriori} error indicator for discontinuous
              {G}alerkin approximations of fourth-order elliptic problems},
   JOURNAL = {IMA J. Numer. Anal.},
  FJOURNAL = {IMA Journal of Numerical Analysis},
    VOLUME = {31},
      YEAR = {2011},
    NUMBER = {1},
     PAGES = {281--298},
      ISSN = {0272-4979},
   MRCLASS = {65N30 (65N15)},
  MRNUMBER = {2755946},
MRREVIEWER = {Manfred Dobrowolski},
       DOI = {10.1093/imanum/drp023},
       URL = {https://doi.org/10.1093/imanum/drp023},
}

@article {CarstensenPuttkammer2023,
    AUTHOR = {Carstensen, Carsten and Puttkammer, Sophie},
     TITLE = {Direct guaranteed lower eigenvalue bounds with optimal a priori convergence rates for the bi-Laplacian},
   JOURNAL = {arXiv},
      YEAR = {2023},
     PAGES = {1--67},
}

@book {ErnGuermond2021,
    AUTHOR = {Ern, Alexandre and Guermond, Jean-Luc},
     TITLE = {Finite elements {I}---{A}pproximation and interpolation},
    SERIES = {Texts in Applied Mathematics},
    VOLUME = {72},
 PUBLISHER = {Springer, Cham},
      YEAR = {2021},
     PAGES = {xii+325},
      ISBN = {978-3-030-56340-0; 978-3-030-56341-7},
   MRCLASS = {65-01},
  MRNUMBER = {4242224},
       DOI = {10.1007/978-3-030-56341-7},
       URL = {https://doi.org/10.1007/978-3-030-56341-7},
}

@ARTICLE{Gallistl2015,
  author = {D. Gallistl},
  title = {{M}orley finite element method for the eigenvalues of the biharmonic
 operator},
  journal = {IMA J. Numer. Anal.},
  year = {2015},
  volume = {35},
  number = {4},
  pages = {1779--1811},
  doi = {10.1093/imanum/dru054},
  eprint = {1406.2876},
  fjournal = {IMA Journal of Numerical Analysis},
  full_text = {https://academic.oup.com/imajna/article-pdf/35/4/1779/5097743/dru054.pdf?guestAccessKey=7ce1b625-d3b0-415d-babb-2f37d806ca08},
}

@article {QiuZhang2020,
    AUTHOR = {Qiu, Weifeng and Zhang, Shun},
     TITLE = {Adaptive first-order system least-squares finite element
              methods for second-order elliptic equations in nondivergence
              form},
   JOURNAL = {SIAM J. Numer. Anal.},
  FJOURNAL = {SIAM Journal on Numerical Analysis},
    VOLUME = {58},
      YEAR = {2020},
    NUMBER = {6},
     PAGES = {3286--3308},
      ISSN = {0036-1429,1095-7170},
   MRCLASS = {65N30 (65N12 65N15 65N50)},
  MRNUMBER = {4173220},
MRREVIEWER = {Liang\ Ge},
       DOI = {10.1137/19M1271099},
       URL = {https://doi.org/10.1137/19M1271099},
}

@article {GuzmanLischkeNeilan2022,
    AUTHOR = {Guzm\'{a}n, Johnny and Lischke, Anna and Neilan, Michael},
     TITLE = {Exact sequences on {W}orsey-{F}arin splits},
   JOURNAL = {Math. Comp.},
  FJOURNAL = {Mathematics of Computation},
    VOLUME = {91},
      YEAR = {2022},
    NUMBER = {338},
     PAGES = {2571--2608},
      ISSN = {0025-5718,1088-6842},
   MRCLASS = {65N30 (65L60 76M10 78M10)},
  MRNUMBER = {4473097},
       DOI = {10.1090/mcom/3746},
       URL = {https://doi.org/10.1090/mcom/3746},
}

@article{g2017b ,
  author = {D. Gallistl},
  title = {Variational formulation and numerical analysis of linear
    elliptic equations in nondivergence form with {C}ordes coefficients},
  journal = {SIAM J. Numer. Anal.},
  year = {2017},
  volume = {55},
  number = {2},
  pages = {737--757},
  doi = {10.1137/16M1080495},
}

@article {DebrabantJakonsen2013,
    AUTHOR = {Debrabant, Kristian and Jakobsen, Espen R.},
     TITLE = {Semi-{L}agrangian schemes for linear and fully non-linear
              diffusion equations},
   JOURNAL = {Math. Comp.},
  FJOURNAL = {Mathematics of Computation},
    VOLUME = {82},
      YEAR = {2013},
    NUMBER = {283},
     PAGES = {1433--1462},
      ISSN = {0025-5718,1088-6842},
   MRCLASS = {65M06 (35K55 49L20)},
  MRNUMBER = {3042570},
MRREVIEWER = {Giacomo\ Dimarco},
       DOI = {10.1090/S0025-5718-2012-02632-9},
       URL = {https://doi.org/10.1090/S0025-5718-2012-02632-9},
}

@article {KaweckiSmears2021,
    AUTHOR = {Kawecki, Ellya L. and Smears, Iain},
     TITLE = {Unified analysis of discontinuous {G}alerkin and
              {$C^0$}-interior penalty finite element methods for
              {H}amilton-{J}acobi-{B}ellman and {I}saacs equations},
   JOURNAL = {ESAIM Math. Model. Numer. Anal.},
  FJOURNAL = {ESAIM. Mathematical Modelling and Numerical Analysis},
    VOLUME = {55},
      YEAR = {2021},
    NUMBER = {2},
     PAGES = {449--478},
      ISSN = {2822-7840,2804-7214},
   MRCLASS = {65N30 (65N12 65N15)},
  MRNUMBER = {4229194},
MRREVIEWER = {Marco\ Verani},
       DOI = {10.1051/m2an/2020081},
       URL = {https://doi.org/10.1051/m2an/2020081},
}

\end{document}